\documentclass[ejs]{imsart}

\RequirePackage[OT1]{fontenc}
\RequirePackage{amsthm,amsmath,amssymb}
\RequirePackage[numbers]{natbib}
\RequirePackage[colorlinks,citecolor=blue,urlcolor=blue]{hyperref}

\pubyear{2006}
\volume{0}
\issue{0}
\firstpage{1}
\lastpage{8}

\startlocaldefs
\numberwithin{equation}{section}
\theoremstyle{plain}
\newtheorem{theorem}{Theorem}[section]
\newtheorem{lemma}{Lemma}[section]
\endlocaldefs

\usepackage{graphicx}

\DeclareMathOperator*{\argmin}{\arg\!\min}

\newcommand{\ms}{\boldsymbol}
\newcommand{\mb}{\mathbf}

\begin{document}

\begin{frontmatter}
\title{Variable Selection via Adaptive False Negative Control in Linear Regression}
\runtitle{False Negative Control in Regression}

\begin{aug}
\author{\fnms{X. Jessie} \snm{Jeng}\thanksref{t1}\ead[label=e1]{xjjeng@ncsu.edu}}
\and
\author{\fnms{Xiongzhi} \snm{Chen}\ead[label=e2]{xiongzhi.chen@wsu.edu}}

\address{Department of Statistics, North Carolina State University, Raleigh, NC 27695\\
\printead{e1}}

\address{Department of Mathematics and Statistics\\ Washington State University, Pullman, WA 99164\\
\printead{e2}\\}

\thankstext{t1}{The research of X. J. Jeng was supported in part by NSF Grant DMS-1811360.}
\runauthor{Jeng and Chen}

\affiliation{Some University and Another University}

\end{aug}

\begin{abstract}
Variable selection methods
have been developed in linear regression to provide sparse solutions. Recent studies have focused on further interpretations on the sparse solutions in terms of false positive control. In this paper, we consider false
negative control for variable selection with the goal to efficiently
select a high proportion of relevant predictors. Different from existing studies in power analysis and sure screening, we propose to directly estimate the  false negative proportion (FNP) of a decision rule and select the smallest subset of predictors that has the estimated FNP less than a user-specified control level.
The proposed method is adaptive to the user-specified control level on FNP by selecting less candidates if a higher level is implemented. On the other hand, when data has stronger effect size or larger sample size, the proposed method controls FNP more efficiently with less false positives.
New analytic techniques are developed to cope with the major challenge of FNP control when relevant predictors cannot be consistently separated from irrelevant ones. Our numerical results are in line with the theoretical findings.
\end{abstract}


\tableofcontents
\end{frontmatter}

\section{Introduction}

We consider a sparse linear model

\begin{equation}
\mathbf{y}=\mathbf{X}\boldsymbol{\beta }+\boldsymbol{\varepsilon },
\label{def:model}
\end{equation}%
where $\mathbf{y}=(y_{1},\ldots ,y_{n})^{T}$ is the vector of $n$
observations of response, $\mathbf{X}=[\mathbf{x}_{1},\ldots ,\mathbf{x}_{p}]\in \mathbb{R}^{n\times p}$ is the design matrix, $\boldsymbol{\beta }%
=(\beta _{1},\ldots ,\beta _{p})^{T}$ is the vector of unknown coefficients,
and $\boldsymbol{\varepsilon }\sim \mathcal{N}_{n}\left( 0,\sigma ^{2}%
\mathbf{I}\right) $ is the vector of random errors.
We assume $\sigma ^{2}=O\left( 1\right) $. Let $I_{1}=\left\{ 1\leq j\leq
p:\beta _{j}\neq 0\right\} $ be the set of indices for non-zero coefficients
with cardinality $s=|I_{1}|$ and $I_{0}=\left\{ 1\leq j\leq p:\beta
_{j}=0\right\} $ with cardinality $p_{0}=\left\vert I_{0}\right\vert$.

Variable selection methods often provide sparse
solutions for the estimation of $\boldsymbol{\beta}$. The non-zero elements
of an estimate correspond to variables selected as candidates for relevant
predictors. A great amount of literature with many fruitful ideas has contributed to the development of sparse solutions to accommodate the underlying features of the data. We refer to \cite{BVdGeer11} and the references therein for a nice introduction.

Given a selection result, a false positive (FP) occurs when an irrelevant
predictor is selected, and a false negative (FN) occurs when a relevant
predictor is not selected. It is natural to interpret a selection result in
terms of false positive or false negative control, and exciting progress has
emerged for false positive control, e.g. \cite{Barber15}, \cite{Bogdan15}, \cite{candes2018panning}, \cite{Gsell16}, \cite{JiZhao14}, \cite{Su16}, \cite{ZhangCheng2017}. However, the
study for efficient false negative control remains relatively underdeveloped.

\textcolor{black}{False negative control is important in many real applications and sometimes a more serious concern than false positive control. For example, in pre-surgical brain mapping with functional MRI, the primary goal is to reduce false negatives where genuine functional areas are not identified.  This is because neurosurgical patients are more likely to experience significant harm from
mistakenly deeming a region to be functionally uninvolved and subsequently resecting critical tissue than from incorrectly assigning
function to an uninvolved region \citep{liu2016pre, loring2002now, silva2018challenges}. Another example where false negative control is of main concern is in the exploratory stage of high-dimensional data analysis, where pre-screening is often conducted to reduce data dimension while keeping a high proportion of true signal variables for follow-up studies.}

The problem of false negative control is conceptually related but
methodologically very different from Sure Screening in, e.g., \cite{fan2008sure, fan2010selective}. Sure Screening
aims to reduce the data dimension by removing
only irrelevant predictors. For instance, the Sure Independence Screening
procedure in \cite{fan2008sure} ranks variables by estimated marginal
regression coefficients and selects the top $d$ variables where $d$ is fixed
at $n-1$ or $n/\log n$. It has been proved that under certain conditions,
the screening procedure has eliminated only irrelevant predictors with high
probability. The false negative control problem considered here
focuses on selecting a high proportion of relevant predictors without
including many unnecessary irrelevant predictors. It may be regarded as a
more refined screening procedure with a data-adaptive selection rule instead
of a fixed $d$.

We use false negative proportion (FNP) as a measure for false negative
control. For a given selection rule, FNP is defined as the ratio of the
number of false negatives to the total number of relevant predictors. FNP
takes values in $[0, 1]$ and is equivalent to $1-\mathsf{Sensitivity}$ in
binary classification framework. Our work starts with consistently
estimating FNP for a given selection rule. To achieve this, we develop novel
analyses on the tail behavior of the empirical processes associated with
FNP.
Based on the estimation of FNP, we develop a new variable selection
procedure to control FNP at a user-specified level. If users can tolerate
more false negatives, they may implement lower control levels on FNP in the
procedure and select less candidates for relevant predictors. On the other
hand, if the effect of relevant predictors gets stronger or sample size
increases, the procedure controls FNP more efficiently with less false positives.

An important component of the proposed FNP control method is an
estimator for the number of relevant predictors. We provide a consistent
estimator for dependent test statistics, for which we adopt the recently
developed debiased Lasso estimates \citep{JM14,vandegeer2014,Zhang:2014}.

\textcolor{black}{Although FNP, by definition, is equivalent to the power in (single) hypothesis testing, our proposed study on FNP control is very different from the existing power analysis in hypothesis testing. In the latter, a decision rule is  built upon Type I error control and followed by power calculation with assumptions on the effect size. For such methods to control FNP in addition to controlling family-wise Type I error when multiple hypotheses are considered, the effect sizes of relevant variables need to be larger enough to ensure essentially perfect separation of relevant and irrelevant variables.
	The proposed method, on the other hand, directly bound the estimated FNP at a user-specified level, which allows a more effective control on FNP. Our condition on effect size for FNP control is shown to be weaker than the existing beta-min conditions that are required for perfect separation of relevant and irrelevant variables.}

The rest of the paper is organized as follows. \autoref{sec:FNP-est}
presents FNP estimation in two steps: (1) constructing test statistics for
regression coefficients and (2) approximating FNP based on the test statistics. \autoref{sec:FNP-control} develops a
variable selection method to control FNP at a user-specified level and a computational algorithm to implement the method.
\autoref{sec:Simulation} presents the finite-sample performance of the proposed method in simulation. Conclusion and further discussion are provided in \autoref{sec:conclude}.
Proofs for the main theoretical results
are presented in \autoref{sec:Proofs}. Extra technical details are
provided in Appendix.

\section{False Negative Proportion Estimation}

\label{sec:FNP-est}

Recall that for a selection rule, FNP is the ratio of the number of false
negatives to the total number of relevant predictors. In this section, we
rank the predictors based on the debiased Lasso estimates and
approximate FNP at a given cut-off point on the list of ranked predictors.

\subsection{Test Statistics Based on Debiased Lasso Estimates}

\label{sec:introDelasso}

Recall model (\ref{def:model}). The well-known Lasso estimator is
\begin{equation}
\hat{\boldsymbol{\beta }}=\hat{\boldsymbol{\beta }}(\lambda )=\arg \min_{%
	\boldsymbol{\beta }\in \mathbb{R}^{p}}(\Vert \mathbf{y}-\mathbf{X}%
\boldsymbol{\beta }\Vert _{2}^{2}/n+2\lambda \Vert \boldsymbol{\beta }\Vert
_{1}),  \label{eqLasso}
\end{equation}%
where $\lambda $ is a tuning parameter \citep{Tibshirani:1996}. Recently, the
debiased Lasso estimator has been developed
to mitigate the bias of Lasso estimator \citep{vandegeer2014, Zhang:2014}.
The debiased Lasso estimator is defined as
\begin{equation}
\hat{\mathbf{b}}=\left( \hat{b}_{1},\ldots ,\hat{b}_{p}\right) ^{T}=\hat{%
	\boldsymbol{\beta }}+n^{-1}\hat{\boldsymbol{\Theta }}\mathbf{X}^{T}(\mathbf{y%
}-\mathbf{X}\hat{\boldsymbol{\beta }}),  \label{eqRelaxedEst}
\end{equation}%
where $\hat{\boldsymbol{\Theta }}\in \mathbb{R}^{p\times p}$ is an estimate
for the precision matrix of the predictors and can be obtained via nodewise
regression on $\mathbf{X}$ as in \cite{meinshausen2006}. Let $\hat{%
	\boldsymbol{\Sigma }}=n^{-1}\mathbf{X}^{T}\mathbf{X}$. It has been shown
that
\begin{equation}
\sqrt{n}(\hat{\mathbf{b}}-\boldsymbol{\beta })=n^{-1/2}\hat{\boldsymbol{%
		\Theta }}\mathbf{X}^{T}\boldsymbol{\varepsilon }-\boldsymbol{\delta }=%
\mathbf{w}-\boldsymbol{\delta },  \label{eq:b_decomp}
\end{equation}%
where
\begin{equation*}
\mathbf{w}|\mathbf{X}\sim \mathcal{N}_{p}(0,\sigma ^{2}\hat{\boldsymbol{%
		\Omega }}),\qquad \hat{\boldsymbol{\Omega }}=\hat{\boldsymbol{\Theta }}\hat{%
	\boldsymbol{\Sigma }}\hat{\boldsymbol{\Theta }}^{T},
\end{equation*}%
and
\begin{equation}
\boldsymbol{\delta }=\left( \delta _{1},\ldots ,\delta _{p}\right) ^{T}=%
\sqrt{n}(\hat{\boldsymbol{\Theta }}\hat{\boldsymbol{\Sigma }}-\mathbf{I})(%
\hat{\boldsymbol{\beta }}-\boldsymbol{\beta }).  \label{eqDelta}
\end{equation}%
Under certain conditions, $\Vert \boldsymbol{\delta }\Vert _{\infty
}=o_{p}(1)$, which implies the asymptotic normality of $\hat{\mathbf{b}}$ %
\citep{CaiGuo17,JM14,  JM2018, vandegeer2014, Zhang:2014}. \textcolor{black}{We present the set of conditions from \cite{JM2018} as A1) - A3) in Appendix \ref{sec:extraDLasso}. }

In this paper, we obtain test statistics for $\boldsymbol{\beta }$ using the
standardized debiased Lasso estimator as
\begin{equation}
z_{j}=\sqrt{n}\hat{b}_{j}\sigma ^{-1}\hat{\boldsymbol{\Omega }}_{jj}^{-1/2}%
\text{ \ for}\qquad 1\leq j\leq p  \label{def:Z}
\end{equation}%
where $\hat{\boldsymbol{\Omega}}_{jj}$ denotes the $\left( j,j\right) $
entry of $\hat{\boldsymbol{\Omega}}$. Therefore, for each $1\leq j\leq p$,
\[
z_{j}=\mu _{j}+w_{j}^{\prime }-\delta _{j}^{\prime },
\]
where, given $\mathbf{X}$,
\begin{equation} \label{defWprime}
w_{j}^{\prime }=\frac{w_{j}}{\sigma \sqrt{\boldsymbol{\hat{\Omega}}_{jj}}}%
\sim \mathcal{N}(0,1),\text{ \ }\delta _{j}^{\prime }=\frac{\delta _{j}}{%
	\sigma \sqrt{\boldsymbol{\hat{\Omega}}_{jj}}}\text{ \ and\ \ \ }\mu _{j}=%
\frac{\sqrt{n}\beta _{j}}{\sigma \sqrt{\boldsymbol{\hat{\Omega}}_{jj}}}.
\end{equation}%


\subsection{Approximating False Negative Proportion}


We aim to determine a cut-off value for the realized test statistics to control false negative proportion (FNP) at an user-specified level. For this purpose, we first study the consistent estimation of $\mathsf{FNP}$.
For any $t>0$, define
\begin{equation*}
R\left( t\right) =\sum_{j=1}^{p}1_{\left\{ \left\vert z_{j}\right\vert
	>t\right\} }, \quad \mathsf{TP}\left( t\right) =\sum_{j\in I_{1}}1_{\left\{
	\left\vert z_{j}\right\vert >t\right\} },  \quad
\mathsf{FN}\left( t\right)
=\sum_{j\in I_{1}}1_{\left\{ \left\vert z_{j}\right\vert \leq t\right\} },
\quad \mathsf{FP}\left( t\right) =\sum_{j\in I_{0}}1_{\left\{ \left\vert
	z_{j}\right\vert >t\right\} }.
\end{equation*}
Note that $\mathsf{FN}\left( t\right) $ is unobservable as $I_{1}$ is
unknown, and that the dependence among $z_{j}$'s also affect $\mathsf{FN}\left(
t\right) $. It is easy to see that
\begin{equation}
\mathsf{FN}\left( t\right) =s-\mathsf{TP}\left( t\right) =s-\left[ R\left(
t\right) -\mathsf{FP}\left( t\right) \right]  \label{eq:TP}
\end{equation}%
and
\begin{equation}
\mathsf{FNP}\left( t\right) =\frac{\mathsf{FN}\left( t\right) }{s}=1-\frac{%
	R\left( t\right) -\mathsf{FP}\left( t\right) }{s}.  \label{eq:FNPa}
\end{equation}
Since $R\left( t\right) $ is directly observable from the data, the unknown
quantities in (\ref{eq:FNPa}) are $\mathsf{FP}\left( t\right) $ and $s$. We
propose to substitute $\mathsf{FP}\left( t\right) $ in (\ref{eq:TP}) by $%
2(p-s)\Phi \left( -t\right) $, where $\Phi (\cdot )$ is the cumulative distribution function (CDF) of a standard Normal random variable,
because $z_{j}$ is asymptotically standard Normal for $j\in I_{0}$. Further, we can plug in
an estimator $\hat{s}$ for $s$, which results in the estimator
\begin{equation}
\widehat{\mathsf{FNP}}\left( t\right) =1-\frac{R\left( t\right) -2(p-\hat s) \Phi
	\left( -t\right) }{\hat{s}}.  \label{def:hatFNP}
\end{equation}

From the definitions of $\mathsf{FNP}\left( t\right) $ and $\widehat{\mathsf{FNP}}\left( t\right) $, it can be shown that \\
$\left\vert \widehat{\mathsf{FNP}}\left(t\right) -\mathsf{FNP}\left(t\right) \right\vert =o_{P}(1)$ is
implied by
\begin{equation}  \label{eq:FPandS}
s^{-1}\left\vert \mathsf{FP}(t)-2p_0\Phi (-t)\right\vert =o_{P}\left( 1\right)
\qquad \text{and} \qquad |\hat s /s - 1| = o_P(1).
\end{equation}
Because $\mathsf{FP}(t)$ is the summation of $p_0$ terms and $s$ can be much smaller than $p_0$, approximating  $\mathsf{FP}(t)/s$ requires more delicate analysis than approximating $\mathsf{FP}(t)/p_0$, which has been studied in the literature for False Discovery Proportion (FDP) control (e.g. \cite{fan2012estimating}). Also, the dependence among test statistics $\{z_j\}_{j=1}^p$ adds another layer of difficulty.

In this paper, we consider $s=p^{1-\eta }$ for some $\eta \in (0, 1)$, so
that the number of relevant predictors is of a smaller order than the total
number of variables. On the other hand, we consider $t$ values calibrated as
$t = t_\xi =\sqrt{2 \xi \log p}$ for some $\xi >0$, so that the scale of $t$
is comparable to that of the extreme value of $p$ independent standard
Gaussian variables. Such calibration has been utilized to study the detection of Gaussian mixtures
\citep{Arias2018, tony2011optimal, donoho2004higher}, and to analyze variable selection consistency in linear regression \citep{ji2012ups}.
In this paper, we adopt the calibration to study the estimation of FNP in linear regression.

Further,  define the precision matrix of the predictors as $\boldsymbol{\Theta }$
and let
\begin{equation*}
s_{j}=\left\vert \left\{ k\neq j:\boldsymbol{\Theta }_{jk}\neq 0\right\}
\right\vert \qquad \text{and} \qquad s_{max }=\max_{1\leq j\leq p}s_{j}.
\end{equation*}
Namely, the parameter $s_{max}$ represents the row-sparsity of the precision matrix, which contributes to the strength of the dependence among the test statistics.
Define
\[
\gamma^*_1 = 2\eta - \min\{1, {\log(n/s_{max}) \over 2 \log p}\}, \qquad \gamma^*_2 = 2 -2\eta -{\log n \over 2 \log p},
\]	
and
\begin{equation} \label{def:gamma*}
\gamma^*  =  \max\{\gamma^*_1, \gamma^*_2\}.
\end{equation}	

The next theoretical result demonstrates the range of $t$ values in which the first equation
$s^{-1}\left\vert \mathsf{FP}(t)-2p\Phi (-t)\right\vert=o_{P}\left( 1\right)$
in (\ref{eq:FPandS}) is achievable.

\begin{theorem}
	\label{thm:FP} Consider model (\ref{def:model}) and the test statistics $\{z_{j}\}_{j=1}^{p}$ in (\ref{def:Z}). Assume conditions A1)
	through A3) in Appendix \ref{sec:extraDLasso} for the asymptotic normality
	of $\{z_{j}\}_{j=1}^{p}$. Let $s=p^{1-\eta }$ for some $\eta \in (0, 1)$ and
	$t=t_{\xi }=\sqrt{2\xi \log p}$ for $\xi >0$. Assume $\xi > \min\{\eta, \gamma^*\}$ for $\gamma^*$ in (\ref{def:gamma*}),
	then
	\begin{equation}  \label{eq:FP(xi)}
	s^{-1}\left\vert \mathsf{FP}(t_{\xi })-2p_0\Phi (-t_{\xi })\right\vert
	=o_{P}\left( 1\right).
	\end{equation}
\end{theorem}	
Because $\mathsf{FP}(t)$ is the summation of $p_0$ indicator functions and $%
p_0 \gg s (=p^{1-\eta})$, $\mathsf{FP}(t)/ s$ blows up at constant $t$. \autoref{thm:FP} says that the approximation of $\mathsf{FP}(t)/ s$ is achievable
for $t$ at the scale of $t_\xi$. This is substantially different from the existing study of FDP control, where the approximation of $\mathsf{FP}(t)/ p_0$ and $R(t)/p_0$ are studied at constant $t$.

The condition $\xi > \min\{\eta, \gamma^*\}$ can be decomposed  as follows. When $\eta \le \gamma^*$, we have $\xi > \eta$, and the claim in (\ref{eq:FP(xi)}) follows by showing  that $s^{-1}\mathsf{FP}(t_{\xi }) = o_p(1) = s^{-1} p_0\Phi (-t_{\xi })$. On the other hand, when $\eta > \gamma^*$ and $\gamma^* < \xi \le \eta$,  more delicate analysis is needed to study the variability of $\mathsf{FP}(t_\xi)$. The condition $\xi > \gamma^*_1$ essentially controls the variability of  $s^{-1}\mathsf{FP_{w'}}(t_\xi)$, where $\bf{w}'$ is the Gaussian component of $\mb{z}$ as in (\ref{defWprime}) and $\mathsf{FP_{w'}}(t_\xi) = \sum_{j\in I_{0}}1_{\left\{ \left\vert w'_{j}\right\vert >t_\xi\right\} }$. The condition $\xi > \gamma^*_2$ controls the cumulative errors caused by the component $\ms{\delta}'$ of $\mb{z}$.

Existing study in  \cite{ji2012ups} has shown optimal phase diagram in $(\xi, \eta)$ for high-dimensional variable selection. Their work, however, focuses on scenarios with $\xi > \eta$. We extend the analysis to the more challenging case with $\gamma^* < \xi \le \eta$, for which we study the variability of  $s^{-1}\mathsf{FP_{w'}}(t_\xi)$ under the dependence of test statistics. Recall the covariance matrix $\sigma^2 \hat{\ms{\Omega}}$ in (\ref{eq:b_decomp}). Since  $\hat{\boldsymbol{\Omega }}=\hat{\boldsymbol{\Theta }}\hat{\boldsymbol{\Sigma }}\hat{\boldsymbol{\Theta }}^{T}$ and that $\hat{\ms{\Sigma}}$ is not a sparse matrix, $\sigma^{2}\hat{\boldsymbol{\Omega}}$ is not sparse or possessing any well-known structures.
The study in \cite{ji2012ups} imposes  conditions on the covariance matrix of predictors that essentially prohibit excessive signal cancellations when performing marginal regression. Our condition of dependence, on the other hand, demonstrate the effect of the sparsity of precision matrix ($s_{max}$) through $\gamma^*$.  Overall, $\xi > \min\{\eta, \gamma^*\}$ is easier to be satisfied with larger $n$, smaller $p$, or smaller $s_{max}$. 	


To achieve the second equation in (\ref{eq:FPandS}), we modify the estimator
introduced in \cite{MR06} and study its consistency for estimating $s$ in our setting. We refer to the modified estimator as the MR estimator.
Recall the standardized debiased Lasso estimator $z_{j}=\sqrt{n}\hat{b}_{j}\sigma
^{-1}\hat{\boldsymbol{\Omega }}_{jj}^{-1/2},1\leq j\leq p$. Let $%
F_{p}(t)=p^{-1}\sum_{j=1}^{p}1_{\{|z_{j}|>t\}}$ and $\bar{\sigma}\left(
t\right) =\sqrt{2\bar{\Phi}(t)\left( 1-2\bar{\Phi}(t)\right) }$, where $\bar{%
	\Phi}(t)=1-\Phi (t)$. The MR estimator for the portion of relevant
predictors ($\pi =s/p$) is constructed as
\begin{equation}
\hat{\pi}=\sup_{t>0}{\frac{F_{p}(t)-2\bar{\Phi}(t)-c_{p}\bar{\sigma}\left(
		t\right) }{1-2\bar{\Phi}(t)}},  \label{def:hat_pi}
\end{equation}%
where $c_{p}$ is a bounding sequence pre-specified as follows. Define $G_{p}\left(  t\right)  =p^{-1} \sum_{j=1}^{p}1_{\{|w_{j}^{\prime}|>t\}}$,
\begin{equation} \label{def:HC}
\mathsf{H}(t)={\frac{G_p(t)-2\bar{%
			\Phi}(t)}{\bar{\sigma}\left( t\right) }}, \qquad \text{and}\qquad
V_{p}=\sup_{t>0}\mathsf{H}(t).
\end{equation}%
Set $c_{p}$ as the $(1-\alpha _{p})$-th quantile of $V_{p}$ for $\alpha_{p}=o(1)$, so that $P(V_{p}>c_{p})=\alpha _{p}\rightarrow 0$ as $%
p\rightarrow \infty $.
In other words, $c_{p}$ can be looked upon as an
upper bound for $V_{p}$ probabilistically, and the implement of $c_p$ in (\ref{def:hat_pi}) eventually controls over-estimation on  $\pi$.

Compared to the original MR estimator in \cite{MR06}, the key
modification in (\ref{def:hat_pi}) and (\ref{def:HC}) is the use of
$F_{p}(t) $ and $G_{p}\left(  t\right)$, two empirical processes each
with dependent random summands. Naturally, this requires different techniques to
find $\left\{  c_{p}\right\}  _{p\geq1}$. The setting in
\cite{MR06} considers independent $p$-values that are uniformly
distributed under the null hypothesis. Since the limiting distribution of the
uniform empirical process with independent summands is known and has an
analytic expression, a bounding sequence can be directly found from the
distribution in the construction of the original MR estimator. However, in our
settings $\left\{  z_{j}\right\}  _{j=1}^{p}$ are dependent, and the exact
distributions of $\{\hat b_{j}\}_{j=1}^{p}$ are unspecified.

In fact, $\{\hat b_{j}\}_{j=1}^{p}$ asymptotically has covariance matrix $\sigma^{2}\hat{\boldsymbol{\Omega}} =\sigma^2\hat{\boldsymbol{\Theta }}\hat{\boldsymbol{\Sigma }}\hat{\boldsymbol{\Theta }}^{T}$. In theory, $|\hat{\boldsymbol{\Omega}}_{ij} - \boldsymbol{\Theta}_{ij}| = o_p(1)$ for any $(i, j)$ under  conditions A1) through A3) in Appendix \ref{sec:extraDLasso}. However, $\hat{\boldsymbol{\Omega}}$ itself is neither diagonal nor sparse, and the approximation errors of all the elements in $\hat{\boldsymbol{\Omega}}$ add up to influence $V_p$ .
Note that $V_p$ is the higher criticism statistic of \cite{donoho2004higher} based on the Gaussian component $\mb{w}'$ of $\mb{z}$.
Unfortunately,  existing techniques for higher criticism statistic under short-range and long-range dependence \citep{hall2008HC} cannot be applied here because our test statistics with covariance matrix $\sigma^{2}\hat{\boldsymbol{\Omega}}$ cannot be  partitioned as in \cite{hall2008HC}.

In this paper, we employ a discretization technique adopted from \cite{arias2011} to derive bounds on the variance of a discretized
$\left\{  \mathsf{H}(t)   :t>0\right\}$ and define
a discretized version of $V_{p}$ as
\begin{equation}
V_{p}^{\ast }=\max \left\{ \mathsf{H}(t)  :t\in \left[ \sqrt{
	\tau _{0}\log p},\sqrt{\tau _{1}\log p}\right] \cap \mathbb{N}\right\}
\label{def:Vp*}
\end{equation}%
for two positive constants $\tau _{0}$ and $\tau _{1}$ such that $0
<\tau_{0}<\tau _{1}$. Then, a discretized version of the MR estimator is defined as
\begin{equation}
\hat{\pi}^{\ast }=\max \left\{ \frac{F_{p}(t)-2\bar{\Phi}(t)-c_{p}^{\ast }%
	\bar{\sigma}\left( t\right) }{1-2\bar{\Phi}(t)}:t\in \left[ \sqrt{\tau
	_{0}\log p},\sqrt{\tau _{1}\log p}\right] \cap \mathbb{N}\right\},
\label{def:pi^*}
\end{equation}
where $c_{p}^{\ast}$ is the $(1-\alpha _{p})$-th quantile of $V_{p}^{\ast}$ for $\alpha _{p} = o(1)$.
Let
\begin{equation}
\mu _{min}=\min_{j\in I_{1}}\sqrt{n}\left\vert \beta _{j}\right\vert \sigma
^{-1}\boldsymbol{\Theta}_{jj}^{-1/2}  \label{def:mu_min}
\end{equation}%
as a measure on the minimal effect size of relevant variables. The following theorem demonstrates the consistency of $\hat{\pi}^{\ast }$.
Its proof is presented in \autoref{sec:proof_consistency}.

\begin{theorem}
	\label{thm:proportion} Assume conditions A1) through A3) in Appendix \ref{sec:extraDLasso} for the asymptotic normality
	of $\{z_{j}\}_{j=1}^{p}$. Let $s=p^{1-\eta }$ for some $\eta \in (0,1)$ and $%
	\mu_{min} \ge \sqrt{2(\gamma^*+c) \log p}$ for $\mu_{min}$ and $\gamma^* $
	in (\ref{def:mu_min}) and (\ref{def:gamma*}) and some constant $c>0 $. Then $%
	\hat{\pi}^{\ast }$ \textcolor{black}{with bounding sequence $c_{p}^{\ast}$ at the order of $(s_{\max}/n)^{1/4} \log p$}, $\tau_0 \in (2\gamma^*, 2\gamma^* + c)$, and $\tau_1 > 2(\gamma^*+c)$
	consistently estimates the proportion $\pi $ of relevant predictors, i.e.,
	for any $\delta >0$,
	\begin{equation*}
	P(|\hat{\pi}^{\ast }/\pi - 1| < \delta)\rightarrow 1
	\end{equation*}
	and, consequently, $\hat s = \hat \pi^* p$ satisfies
	\begin{equation*}
	P(|\hat{s}/s -1| <\delta)\rightarrow 1.
	\end{equation*}
\end{theorem}
Note that
the order of $c^*_p$ shows the effects of sparsity of the precision
matrix ($s_{max}$), sample size ($n$), and dimensionality ($p$) on $V_p^*$.
The condition $\mu_{min} \ge \sqrt{2(\gamma^*+c) \log p}$ shows that
consistent estimation of $s$ gets easier with smaller $s_{max}$, larger $n$, and smaller $p$.

In summary, \autoref{thm:FP} and \autoref{thm:proportion} facilitate the two equations in (\ref{eq:FPandS}) for FNP$(t)$ estimation by $\widehat{\mathsf{FNP}}(t)$. \textcolor{black}{Note that in practice we will need to simulate $V_{p}$ and $c_p$ to derive the estimated $s$ and $\mathsf{FNP}(t)$. Please refer to \autoref{sec:Algorithm} for details of the numerical implementation.}


\section{FNP Control at a User-Specified Level}

\label{sec:FNP-control}

In this section, we introduce a new method for FNP control at a user-specified level in high-dimensional regression.
We say that a variable selection method asymptotically controls FNP at a
pre-specified level $\epsilon \in (0, 1) $ if the FNP of its selection outcome
satisfies
\begin{equation*}
P\left( \mathsf{FNP}<\epsilon \right) \rightarrow 1.
\end{equation*}
Such methods are useful in applications where data dimensions need to be
largely reduced for subsequent analyses while controlling false negatives at
a tolerable level.

\subsection{The FNC-Reg Procedure}

Based on the approximation results of FNP, we propose the False Negative Control for Regression (FNC-Reg) procedure,
which determines the cut-off threshold on the list of ranked $\{|z_{j}|\}_{j=1}^{p}$ as
\begin{equation}
t^{\ast }(\epsilon )=\sup \left\{ t:\widehat{\mathsf{FNP}}(t)\leq \epsilon
\right\}  \label{def:t*}
\end{equation}
for an user-specified $\epsilon \in (0, 1)$. FNC-Reg selects predictors with $|z_{j}|>t^{\ast }(\epsilon )$.

\textcolor{black}{It can be seen that FNC-Reg is a procedure built upon direct estimation of FNP and a user-specified control level of FNP. Given that $\widehat{\mathsf{FNP}}(t)$ is non-increasing with $t$, FNC-Reg selects the smallest subset of $\{z_{j}\}_{j=1}^{p}$  such that the estimated FNP is less than $\epsilon$. Moreover, this procedure depends on user's preference for the control level of FNP. Since $t^{\ast }(\epsilon )$ is non-decreasing with $\epsilon $, if users can tolerate missing a higher proportion (larger $\epsilon $) of relevant variables, they may select less variables using the procedure. The selected subset of variables can be much smaller than the full set of variables, which corresponds to no false negatives. The next theorem shows that under certain conditions, the FNC-Reg procedure asymptotically controls the true FNP at the level of $\epsilon$.}

\begin{theorem}
	\label{thm:screening} Assume conditions A1) through A3) in Appendix \ref%
	{sec:extraDLasso}. Assume $\mu_{min} \ge \sqrt{2(\gamma^*+c) \log p}$ for $\mu_{min}$ and $\gamma^*$ in (\ref{def:mu_min}) and (\ref{def:gamma*}) and some constant $c>0$. Then $t^*(\epsilon)$
	determined by (\ref{def:t*}) with $\hat s = \hat \pi^*p $ satisfies
	\begin{equation}
	P\left( \mathsf{FNP}(t^{\ast }(\epsilon )) \le \epsilon \right) \rightarrow 1.
	\label{eq:screen}
	\end{equation}
	Consequently, $t^*(\epsilon)$
	determined by (\ref{def:t*}) with $\hat s = \hat \pi p $ and $c_p = c_p^*$ also satisfies
	\begin{equation} \label{eq:screen1}
	P\left( \mathsf{FNP}(t^{\ast }(\epsilon )) \le \epsilon \right) \rightarrow 1.
	\end{equation}
\end{theorem}
Result in (\ref{eq:screen}) shows the FNP control by FNC-Reg when the discretized MR estimator is implemented in (\ref{def:hatFNP}). (\ref{eq:screen1}) extends the result to FNC-Reg with the MR estimator.



We compare the condition on $\mu_{\min}$ in \autoref{thm:screening} with the beta-min condition of variable selection consistency. Our condition on $\mu_{\min}$ achieves the order $O(\sqrt{(\log p) / n})$ for $\beta_{min}$, which is the optimal oder for variable selection consistency \citep{JM2013, wainwright2009}.
On the other hand, Our condition on $\mu_{\min}$ specifies the constant term $\sqrt{2\gamma^*}$ with $\gamma^*$ in (\ref{def:gamma*}),  while existing beta-min conditions for different methods have various constant terms that are often not fully specified.
Therefore, we attempt to compare with the optimal constant term for variable selection consistency in the ideal setting, where the predictors $(X_{i1}, \ldots, X_{ip})$ are generated as i.i.d. samples from $N(0, I_{p\times p})$. Existing study in, for example, \cite{ji2012ups} has shown that the optimal constant is  $\sqrt{2} + \sqrt{2(1-\eta)}$,  i.e. smaller $\eta$ (larger $s$) makes it harder to perfectly separate all the signals from noise.
Then, it follows that $\sqrt{2\gamma^*} < \sqrt{2} + \sqrt{2(1-\eta)}$ for any $\eta \in (0, 1)$. The above analysis shows that in the ideal setting, our condition on $\mu_{\min}$ is  weaker than the optimal beta-min condition for variable selection consistency, and that FNP control can be achieved by FNC-Reg when relevant and irrelevant variables may not be perfectly separated.

\subsection{Numerical Implementation of FNC-Reg}
\label{sec:Algorithm}

We provide a computational algorithm to implement the proposed FNC-Reg
procedure. First, the estimation of $s$ relies on the bounding sequence $c_p$, which is pre-fixed as the $(1-\alpha_p)$-th quantile of $V_p$. 
\textcolor{black}{In numerical implementation, we suggest to
simulate $V_p$ and $c_p$ as follows. We simulate the data under the global null
hypothesis that no relevant predictors exist and calculate the corresponding  standardized
debiased Lasso estimator $\tilde z_j$.}
Note that $\tilde z_j$ is asymptotically distributed as  $w'_j$ under the global null hypothesis. We order $\tilde z_j$'s by their absolute values such that $%
|\tilde z_{(1)}| > |\tilde z_{(2)}| > \ldots > |\tilde z_{(p)}|$ and calculate
\begin{equation}  \label{def:tilde_V}
\tilde{V_p} = \max_{1<j <p/2} {\frac{j/p - 2 \bar{\Phi}(|\tilde z_{(j)}|) }{\bar{\sigma}\left( |\tilde z_{(j)}|\right)}}.
\end{equation}
Repeat the above 1000 times and determine $\tilde{c}_p$ as the $(1-1/\sqrt{%
	\log p})$-th quantile of the empirical distribution of $\tilde{V}_p^{(1)},
\ldots, \tilde{V}_p^{(1000)}$. Consequently, given the ordered test statistics  $|z_{(1)}| > |z_{(2)}| > \ldots > |z_{(p)}|$,
calculate
\begin{equation}  \label{def:tilde_pi}
\tilde{\pi}=\max_{1<j <p/2}{\frac{j/p - 2\bar{\Phi}(|z_{(j)}|)-\tilde c_{p}%
		\bar{\sigma} \left( |z_{(j)}|\right) }{1-2\bar{\Phi}(|z_{(j)}|)}}.
\end{equation}

\vspace{0.1in}

\noindent \underline{{\bf Algorithm 1} FNC-Reg}

\begin{enumerate}
	\item Derive the debiased Lasso estimator $\hat{\mb{b}}$ as in (\ref{eqRelaxedEst}).
	\item Standardize $\hat{\mb{b}}$ and obtain $z_j = \sqrt{n} \hat b_j \sigma^{-1} \hat{\ms{\Omega}}_{jj}^{-1/2}$ for $1 \le j \le p$. Order the $\{z_j\}_{j=1}^p$ as $|z_{(1)}| > |z_{(2)}| > \ldots > |z_{(p)}|$.
	\item Calculate the bounding sequence $\tilde{c}_p$ as the $(1-1/\sqrt{\log p})$-th quantile of the empirical distribution of $\tilde{V}_p$ in (\ref{def:tilde_V}).
	\item Obtain $\tilde \pi$ by (\ref{def:tilde_pi}) and $\hat s = \tilde \pi p$.
	\item Calculate $\widehat{\mathsf{FNP}}(|z_{(j)}|)$ for $j = 1, \ldots, p$ by (\ref{def:hatFNP}).
	\item Obtain $t^*(\epsilon) = \max\{|z_{(j)}|: \widehat{\mathsf{FNP}}(|z_{(j)}|) \le \epsilon\}$ for a user-specified $\epsilon>0$.
	\item Select predictors with $|z_{j}| \ge t^*(\epsilon)$.
\end{enumerate}

\section{Numerical Analysis}
\label{sec:Simulation}

Examples in this section have the response $y$ simulated by the regression
model (\ref{def:model}) with $\boldsymbol{\varepsilon} \sim \mathcal{N}%
_{n}(0,\mathbf{I})$. Each row of $\mathbf{X}$ is simulated from $\mathcal{N}%
_p(0, \boldsymbol{\Sigma})$. We use the Erg{\"o}s-R{\'e}nyi random graph in
\cite{CLZ16} to generate the precision matrix $\boldsymbol{\Theta} =
\boldsymbol{\Sigma}^{-1}$ with $s_{max} \sim \mathsf{Binomial}(p, \theta)$,
such that the nonzero elements of $\boldsymbol{\Theta}$ are randomly located
in each of its rows with magnitudes randomly generated from the uniform
distribution $\mathsf{Uniform}[0.4, 0.8]$. The nonzero coefficients are set at $\beta_1, \ldots, \beta_s$ with the same values. The
debiased Lasso estimates are obtained by applying the R package \textrm{\textit{hdi}}
\citep{hdi2015}.

\subsection{Estimating $s$}

\begin{figure}[!h]
	\caption{Box-plots of $\hat s/s$ with $p=200, s= 10$, and $\protect\beta_1$ increasing from	$0.2$ to $0.5$. The left plot has $n=100$ and the right plot has $n=150$.}
	\label{fig:FNP}
	\includegraphics[width=0.45\textwidth,
	height=0.32\textheight]{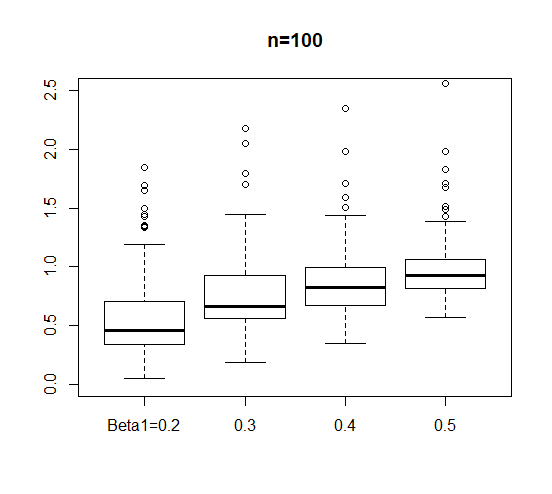}
	\hspace{0.1in}
	\includegraphics[width=0.45\textwidth,
	height=0.32\textheight]{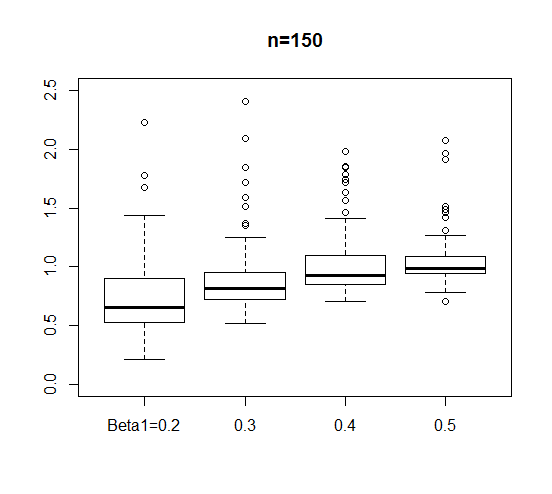}
\end{figure}

We compare the estimated $\hat s$ with the true $s$ in two settings. The
first setting has $p=200,n=100,s=10, \theta=0.02$, and $\beta _{1}= 0.2 - 0.5$. The
second setting increases sample size $n$ to $150$. As claimed in \autoref
{thm:proportion}, the accuracy of $\hat s$ increases with the magnitude of
non-zero coefficients and the sample size. \autoref{fig:FNP} presents the
box-plots of the ratio $\hat s/s$ from 100 replications.
When $\beta_1$ or sample size is small, $\hat s$ tends to under-estimate the true $s$. As $\beta_1$ increases from 0.2 to 0.5 or $n$ increases from 100 to 150,
$\hat s/s$ concentrates more around 1.

\subsection{FNP control}

\label{sec:Simulation_FNPcontrol}

We apply the FNC-Reg algorithm presented in \autoref{sec:Algorithm} to
the simulated data with $p=200$, $n=150$, $s=10$, $\theta = 0.02$, and $\beta_1 = 0.2 - 0.5$. \autoref{tab:FNP_FDP_F} has $\epsilon$ fixed at 0.1 and reports the mean value of $\mathsf{FNP}(t^*)$ as $\beta_1$ increased from 0.2 to 0.5. We also calculated the associated false discovery
proportion ($\mathsf{FDP}(t^{\ast })=\mathsf{FP}(t^{\ast
})/R(t^{\ast })$) to reveal the price in incurring
false positives for FNP control. \textcolor{black}{Further, we calculate the F-measure, which summarizes FNP and FDP by the harmonic mean of (1-FNP) and (1-FDP) \citep{F-measure2011}. F-measure takes a value between 0 and 1, and higher value corresponds to better summarized performance.}

\textcolor{black}{Because we are not aware of any existing methods that directly control FNP in high-dimensional regression, we present the corresponding results of two other methods that perform variable selection based on different criteria. These results help to better understand the results of FNC-Reg. The first method is Lasso whose solution is obtained using the R package {\it hdi}, in which $\lambda$ is determined by cross validation. The second method is Knockoff, which has been developed to control false discovery rate (FDR) at a user-specified level in high-dimensional regression \citep{Barber15, candes2018panning}. We use the "knockoff.filter" function in default from the R package {\it knockoff}, which creates model-X second-order Gaussian knockoffs as introduced in \cite{candes2018panning}. The nominal level is set at 0.1.}

\begin{table}[!h]
	\caption{The mean values and standard deviations (in brackets) of FNP, FDP, and the F-measure from 100 replications for FNC-Reg, Lasso, and Knockoff.}
	\label{tab:FNP_FDP_F}
	\centering
	\vspace{0.1in}
	\begin{tabular}{ll|cccccc}
		\hline
		$\beta_1$ & Method & FNP & FDP & F-measure\\
		[0.05in]
		\hline
		$0.2$ & FNC-Reg & 0.37 (0.16) & 0.35 (0.30) & 0.58 (0.20)  \\
		& Lasso   & 0.29 (0.12) & 0.51 (0.16) & 0.56 (0.16) \\
		& Knockoff & 0.91 (0.23) & 0.06 (0.16) & 0.08 (0.21) \\
		[0.05in]
		\hline
		$0.3$ & FNC-Reg & 0.19 (0.11) & 0.20 (0.24) & 0.77 (0.17)\\
		& Lasso   & 0.15 (0.09) & 0.43 (0.11) & 0.67 (0.09) \\
		& Knockoff & 0.60 (0.44) & 0.11 (0.15) & 0.37 (0.41) \\
		[0.05in]
		\hline
		$0.4$ & FNC-Reg & 0.10 (0.09) & 0.17 (0.25) & 0.84 (0.19)\\
		& Lasso   & 0.06 (0.07) & 0.42 (0.12) & 0.71 (0.10)\\
		& Knockoff & 0.34 (0.45) & 0.12 (0.14) & 0.61 (0.42) \\
		[0.05in]
		\hline
		$0.5$ & FNC-Reg & 0.04 (0.09) & 0.13 (0.22) & 0.90 (0.18)\\
		& Lasso   & 0.02 (0.10) & 0.38 (0.12) & 0.74 (0.14)\\
		& Knockoff & 0.21 (0.49) & 0.09 (0.11) & 0.74 (0.38)\\
		\hline
	\end{tabular}
\end{table}

\textcolor{black}{It can be seen from \autoref{tab:FNP_FDP_F} that as $\beta_1$ increases, the FNP of FNC-Reg decreases, which agrees with the theoretical insight provided by the condition on $\mu_{min}$ in Theorem 3.1. In the challenging scenarios where $\beta_1$ is very small, the FNP of FNC-Reg mostly exceeds the nominal level of 0.1, which is due to the under-estimation of $\hat s$ and FNC-Reg's tendency to select less variables to capture the under-estimated number of signals.
Furthermore, both FNP and FDP of FNC-Reg get smaller for larger $\beta_1$, suggesting that FNC-Reg automatically adapt to and benefit from increasing signal intensity for both false negative and false positive control. \autoref{tab:FNP_FDP_F} also shows that Lasso has lower FNP and much higher FDP than FNC-Reg, which agrees with Lasso's known tendency of over-selection when $p>n$. On the other hand, Knockoff has FDP reasonably  controlled at the nominal level of 0.1 but much higher FNP than those of FNC-Reg and Lasso. In terms of the F-measure that summarizes FNP and FDP, FNC-Reg seems to outperform the other two methods under different   $\beta_1$ values.}

We further illustrate the adaptivity of FNC-Reg to the user-specified control level of FNP. For various values of $\epsilon $, we calculate the relative frequency of the event $\left\{\mathsf{FNP}(t^{\ast }(\epsilon))\le\epsilon \right\}$. \autoref{tab:FNP_epsilon_new} summarizes the results for different settings
with $\epsilon =0.1,0.2,0.3$ and $\beta _{1}=0.3, 0.5, 0.7$. It can be seen
that the relative frequency of $\mathsf{FNP}\le\epsilon $ for FNC-Reg
increases with $\beta_1$, which is consistent with the theoretical insight in \autoref{thm:screening}. On the other hand, for a fixed $\beta_1$, the relative frequency of $\mathsf{FNP}\le\epsilon $ and the
FDP of FNC-Reg decreases with $\epsilon $, which agrees with our expectation for FNC-Reg as more liberal control of FNP incurs less price in false positives. Note that the results of Lasso and Knockoff do not change with the varying $\epsilon$.

\begin{table}[!h]
	\caption{The relative frequencies of the event $\{\mathsf{FNP} < \protect\epsilon\}$ and mean values of FDP from 100 replications for FNC-Reg.}
	\label{tab:FNP_epsilon_new}
	\centering
	\begin{tabular}{l|ccc|ccc|ccc}
		\hline
		& \multicolumn{3}{c|}{$\beta_1= 0.3$} & \multicolumn{3}{c|}{$\beta_1= 0.5$} & \multicolumn{3}{c}{$\beta_1= 0.7$}\\
		& $\epsilon = 0.1$  & $0.2$ & $ 0.3$ & $\epsilon = 0.1$  & $0.2$ & $0.3$ & $\epsilon = 0.1$  & $0.2$ & $0.3$ \\
		[0.05in] \hline
		$1(\mathsf{FNP}\le\epsilon)$ & 0.38 & 0.53 & 0.71 & 0.72 & 0.86 & 0.91 & 0.98 & 0.98 & 0.98 \\
		FDP &  0.20 & 0.15 & 0.09 & 0.13 & 0.08 & 0.05 & 0.14 & 0.10 & 0.07 \\
		\hline
	\end{tabular}%
\end{table}

\section{Conclusion and Discussion}	\label{sec:conclude}

We propose a new variable selection method, FNC-Reg, to efficiently control false negatives in linear regression. Different from existing methods and theory for power analysis and Sure Screening, our procedure directly estimates the FNP of a decision rule and selects the smallest subset of variables that has the estimated FNP less than a user-specified control level.
FNP control is specifically challenging when relevant variables cannot be consistently separated from irrelevant ones due to limited sample size and effect size.
We develop new techniques to analyze FNP control in the challenging setting and to cope with difficulties caused by the dependence of test statistics.

FNC-Reg possesses two types of adaptivity property. First, it adapts to the user's preference level on the control of FNP. When a user can tolerate a less stringent control on FNP,  he or she can input a larger $\epsilon$ in the FNC-Reg procedure and select less variables with less false positives. Secondly, the proposed method is adaptive to the unknown effect size. Note that the implementation of the procedure does not requires the information of effect size. Nevertheless,  the result of the procedure automatically improves in both FNP and FDP as effect size increases.

Our theoretical study presents a weaker condition on $\mu_{min}$ for FNP control by FNC-Reg than the beta-min condition for variable selection consistency. It is also of interest to understand the result of FNC-Reg if the condition on $\mu_{\min}$ may not be satisfied. Assume that among the $s$ signal variables only $s_1$ of them satisfy $\mu_j \ge \sqrt{2 (\gamma^*+c) \log p}$ for some constant $c>0$. Then, similar arguments as in the proof of \autoref{thm:proportion} can be applied to show that $P((1-\delta) s_1 < \hat s < (1+\delta) s) \to 1$ for any $\delta>0$. Note that $\hat s$ does not consistently estimate $s$ anymore,
nor is it a consistent estimator for $s_1$. Such $\hat s$ tends to under-estimate $s$, which can cause the proposed method to select less variables to capture the under-estimated number of signals. Because FNC-Reg ranks the test statistics by their significance and select variables from the top, one can make a statement about FNP control for the signals with effect sizes larger than the observed cut-off position. Such interpretation of results remains valid whether  the condition on $\mu_{min}$ holds or not.

Last but not least, we adopt the debiased Lasso estimator as the test statistic in the paper to demonstrate the new analytic framework of FNP control. We expect that the proposed framework can incorporate other test statistics in linear regression and promote further developments in false negative control based variable selection.

\section{Proofs} \label{sec:Proofs}

This section contains the proofs of \autoref{thm:FP},  \autoref{thm:proportion}, and \autoref{thm:screening}.
Auxiliary lemmas are provided in the appendices. We will frequently use the
Mill's ratio, i.e.,
\begin{equation*}
\bar{\Phi}\left( x\right) ={x}^{-1}\phi \left( x\right) \left( 1+o\left(
1\right) \right) \text{ \ for \ }x\rightarrow +\infty ,
\end{equation*}%
without mentioning it at each instance. All arguments will be conditional on
$\mathbf{X}$, and the symbol $C$ denotes a generic, finite constant whose
values can be different at different occurrences.

\subsection{Proof of \autoref{thm:FP}} \label{sec:proof_thm2.1}

The proof is composed of two parts. The first part assumes $\xi > \eta$ and the second part assumes $\xi\le \eta$.

Consider the first part with $\xi > \eta$.  It suffices to show
$s^{-1}\mathsf{FP}\left( t_{\xi }\right)
=o_{P}(1)$ and $s^{-1}p_0\Phi \left( -t_{\xi }\right) =o(1)$ with $\xi > \eta$. By Mill's ratio,
\begin{equation*}
s^{-1}p\Phi \left( -t_{\xi }\right) \leq Cp^{\eta -\xi }/\sqrt{\log p}=o(1)
\end{equation*}%
when $\xi >\eta $. On the other hand, for a fixed constant $a>0$,
\begin{equation*}
P\left( s^{-1}\mathsf{FP}\left( t_{\xi }\right) >a\right) \leq \frac{E\left(
	\mathsf{FP}\left( t_{\xi }\right) \right) }{as}=\frac{p_{0}\max_{j\in
		I_{0}}P\left( \left\vert z_{j}\right\vert >t_{\xi }\right) }{as}.
\end{equation*}%
The following lemma help quantify the order of $P(|z_{j}|>t_{\xi })$ for $j\in
I_{0}$, and its proof is provided in \autoref{sec:proof_zw-difference}.

\begin{lemma} \label{lemma:zw-difference}
	Assume A1) through A3). Define
	\[
	d_{p}=C_{1}\left( \sqrt{s}\log p/\sqrt{n}+\min \left\{ s,s_{max }\right\}
	\log p/\sqrt{n}\right)
	\] for some constant $C_{1}\geq \max \left\{ 1,2(\sigma \sqrt{C_{\min }}%
	)^{-1}\right\}$. Then
	\[
	\left\vert P(|w'_{j} - \delta'_j|>t_{\xi}) - P(|w'_{j}|>t_{\xi })\right\vert \le C p^{-\xi} d_p + C p^{-2}.
	\]
\end{lemma}
Recall $s=p^{1-\eta} $ with $0<\eta <1$. Then \autoref{lemma:zw-difference}
implies
\begin{equation*}
P\left( s^{-1}\mathsf{FP}\left( t_{\xi }\right) >a\right) \leq \frac{%
	Cp_{0}\left( p^{-\xi} + p^{-\xi} d_p + p^{-2} \right) }{p^{1-\eta}}%
=o\left( 1\right),
\end{equation*}%
where the last step is by $d_p = o(1)$ under A3) and $\xi > \eta$. Then $s^{-1}\mathsf{FP}\left( t_{\xi
}\right) =o_{P}\left( 1\right) $. This justifies the claim of \autoref{thm:FP} for $\xi > \eta$.

Next, we present the second part of the proof with $\xi \le \eta$. Define
\[
D_p = s^{-2}(p^{1-\xi }+p^{2-\xi}\lambda _{1}\sqrt{s_{max}}) \log p.
\]
By the order of $\lambda_1$ in A2) and condition $\xi > \gamma^*_1$, $D_p = o(1)$. Then it is sufficient to show
\[
P(s^{-1}\left\vert \mathsf{FP}(t_{\xi })-2p_0\Phi (-t_{\xi })\right\vert > \sqrt{D_p}) \to 0.
\]
Perform the decomposition
\begin{eqnarray*}
	s^{-1}\left\vert \mathsf{FP}(t_{\xi })-2p_0\Phi (-t_{\xi })\right\vert & \le &  s^{-1}\left\vert \mathsf{FP}(t_{\xi })-\mathsf{E}(\mathsf{FP}(t_\xi))\right\vert \\
	& + & s^{-1}\left\vert \mathsf{E}(\mathsf{FP}(t_\xi) -2p_0\Phi (-t_{\xi }))\right\vert.
\end{eqnarray*}
Then it is sufficient to show
\begin{equation} \label{eq:1}
s^{-1}\left\vert \mathsf{FP}(t_{\xi })-\mathsf{E}(\mathsf{FP}(t_\xi))\right\vert = o_p(\sqrt{D_p})
\end{equation}
and
\begin{equation} \label{eq:2}
s^{-1}\left\vert \mathsf{E}(\mathsf{FP}(t_\xi) -2p_0\Phi (-t_{\xi }))\right\vert = o(\sqrt{D_p}).
\end{equation}

Consider (\ref{eq:1}) first. By Chebyshev's inequality,
\[
P(s^{-1}\left\vert \mathsf{FP}(t_{\xi })-\mathsf{E}(\mathsf{FP}(t_\xi))\right\vert > \sqrt{D_p}) \le {\mathsf{Var}(\mathsf{FP}(t_{\xi })) \over s^2 D_p}.
\]
We derive the order of $\mathsf{Var}(\mathsf{FP}(t_{\xi }))$. By \autoref{lemma:zw-difference},
$P(|w'_{j} - \delta'_j|>t_{\xi}) = P(|w'_{j}|>t_{\xi }) (1+o(1))$ given $d_p = o(1/\sqrt{\log p})$ from A3) and $\xi \le \eta <1$, then  direct calculation gives
\[
\mathsf{Var}(\mathsf{FP}(t_{\xi })) = \mathsf{Var}(\mathsf{FP}_{w'}(t_{\xi }))(1+o(1)),
\]
where $\mathsf{FP}_{w'}(t_{\xi }) = \sum\nolimits_{j \in I_0} 1_{\{|w_{j}^{\prime }|>t_{\xi }\}}$. The following lemma is proved in \autoref{sec:proof_lemma_Var(FP)}.

\begin{lemma}
	\label{lemma:Var(FP)} Assume A1) and A2) and let $t_{\xi }=\sqrt{2\xi \log p}
	$ for any $\xi >0$. Then%
	\begin{equation*}
	\mathsf{Var}\left(
	\sum\nolimits_{j=1}^{p}1_{\{|w_{j}^{\prime }|>t_{\xi }\}}\right)  =O\left( p^{1-\xi }+p^{2-\xi
	}\lambda _{1}\sqrt{s_{max }}\right).
	\end{equation*}
\end{lemma}
The above gives
\[
{\mathsf{Var}(\mathsf{FP}(t_{\xi })) \over s^2 D_p} =  {\mathsf{Var}(\mathsf{FP}_{w'}(t_{\xi })) \over s^2 D_p} (1+o(1)) = o(1),
\]
so that
\[
P(s^{-1}\left\vert \mathsf{FP}(t_{\xi })-\mathsf{E}(\mathsf{FP}(t_\xi))\right\vert > \sqrt{D_p}) \to 0.
\]	

Next consider (\ref{eq:2}). By \autoref{lemma:zw-difference}
\begin{eqnarray*}
	s^{-1}\left\vert \mathsf{E}(\mathsf{FP}(t_\xi) -2p_0\Phi (-t_{\xi }))\right\vert & \le & s^{-1} \sum_{j\in I_0} \left\vert P(|w'_{j} - \delta'_j|>t_{\xi}) - P(|w'_{j}|>t_{\xi })\right\vert \\
	& \le &  C s^{-1} p^{1-\xi} d_p + C s^{-1} p^{-1}.
\end{eqnarray*}
Recall $s = p^{1-\eta}$ and the definitions of $d_p$ and $D_p$. Note that $d_p > s \log p /\sqrt{n}$, then direct calculation gives
\[
s^{-1} p^{1-\xi} d_p = o(\sqrt{D_p})
\]
under condition $\xi>\gamma^*_2$, and
\[
s^{-1} p^{-1} = o(\sqrt{D_p})
\]
with $\xi \le \eta <1$. (\ref{eq:2}) follows consequently.  This concludes the second part of the proof with $\xi \le \eta$.

\subsection{Proof of \autoref{thm:proportion}}
\label{sec:proof_consistency}

First, we have the following lemma showing the order of the bounding sequence $c_{p}^{\ast }$ for $V_{p}^{\ast }$. The proof is provided in \autoref{sec:proof_lemma_boundingSeq}.

\begin{lemma}
	\label{lemma:boundingSeq} Assume conditions A1) through A3) in Appendix \ref{sec:extraDLasso}.
	Consider $V_{p}^{\ast }$ as in (\ref{def:Vp*}). Then $c_{p}^{\ast
	}$ at the order of $(s_{max}/n)^{1/4} \log p$ satisfies $P(V_{p}^{\ast }>c_{p}^{\ast
	})\rightarrow 0$ as $p\rightarrow \infty $.	
\end{lemma}

Now, recall $F_{p}(t)=p^{-1}\sum_{j=1}^{p}1_{\{|z_{j}|>t\}}$ and define $\bar{\Phi%
}_{p}(t)=p^{-1}\sum_{j=1}^{p}1_{\{|\mu _{j}+w_{j}^{\prime }|>t\}}$. Consider
the decomposition
\begin{equation}
\hat{\pi}^{\ast }=\max_{t\in \mathbb{T}}\left\{ {\frac{F_{p}(t)-\bar{\Phi}%
		_{p}(t)}{1-2\bar{\Phi}(t)}}+\frac{\bar{\Phi}_{p}(t)-2\bar{\Phi}%
	(t)-c_{p}^{\ast }\bar{\sigma}\left( t\right) }{1-2\bar{\Phi}(t)}\right\},
\label{eqe13}
\end{equation}%
where $\mathbb{T}$ is defined in (\ref{def:T}).
The first summand within the parentheses on the right hand side (RHS) of (%
\ref{eqe13}) can be safely ignored when bounding $\hat{\pi}^{\ast }/\pi$ 	
as asserted by the following \autoref{lmLast}.

\begin{lemma}
	\label{lmLast}Assume $t = t_{\xi }=\sqrt{2\xi \log p}$ with $\xi > \gamma^* $.
	Then
	\[
	\pi ^{-1} \left\vert F_{p}(t)-\bar{\Phi}_{p}(t)\right\vert \left( 1-2%
	\bar{\Phi}(t)\right) ^{-1}=o_{P}(1).
	\]
\end{lemma}

Define
\begin{equation*}
\hat{\pi}^{\ast \ast }=\max_{t\in \mathbb{T}}\frac{\bar{\Phi}_{p}(t)-2\bar{%
		\Phi}(t)-c_{p}^{\ast }\bar{\sigma}\left( t\right) }{1-2\bar{\Phi}(t)}.
\end{equation*}%
Then it suffices to show
\begin{equation}
P(1-\delta <\hat{\pi}^{\ast \ast }/\pi <1)\rightarrow 1.  \label{9}
\end{equation}

We first show that $\hat{\pi}^{\ast \ast }$ is an asymptotic lower bound of $%
\pi $. Recall the definition of $V_{p}^{\ast }$ as%
\begin{equation*}
V_{p}^{\ast }=\max_{t\in \mathbb{T}}\frac{p^{-1}\sum\nolimits_{j=1}^{p}1_{%
		\{|w_{j}^{\prime }|>t\}}-2\bar{\Phi}(t)}{\bar{\sigma}(t)}.
\end{equation*}%
Since
\begin{equation*}
\bar{\Phi}_{p}(t)\leq p^{-1}s+p^{-1}\sum_{j\in I_{0}}1_{\{|w_{j}^{\prime
	}|>t\}}=\pi +(1-\pi )p_{0}^{-1}\sum_{j\in I_{0}}1_{\{|w_{j}^{\prime }|>t\}},
\end{equation*}%
then%
\begin{eqnarray*}
	P(\hat{\pi}^{\ast \ast }>\pi ) &\leq &P\left( \max_{t\in \mathbb{T}}\left\{
	(1-\pi )\left( p_{0}^{-1}\sum\nolimits_{j\in I_{0}}1_{\{|w_{j}^{\prime
		}|>t\}}-2\bar{\Phi}(t)\right) -c_{p}^{\ast }\bar{\sigma}(t)\right\} >0\right)
	\\
	&\leq &P\left( \max_{t\in \mathbb{T}}\left\{ p_{0}^{-1}\sum\nolimits_{j\in
		I_{0}}1_{\{|w_{j}^{\prime }|>t\}}-2\bar{\Phi}(t)-c_{p_{0}}^{\ast }\bar{\sigma%
	}(t)\right\} >0\right) \\
	&\leq &P\left( V_{p_{0}}^{\ast }>c_{p_{0}}^{\ast }\right) ,
\end{eqnarray*}%
where the second inequality follows since $c_{p}^{\ast }$ is non-decreasing
in $p$ and $c_{p}^{\ast }(1-\pi )^{-1}>c_{p_{0}}^{\ast }$. However,
\autoref{lemma:boundingSeq} asserts $P(V_{p_{0}}^{\ast }>c_{p_{0}}^{\ast
})\rightarrow 0$. So,
\begin{equation}
P(\hat{\pi}^{\ast \ast }>\pi )\leq P(V_{p_{0}}^{\ast }>c_{p_{0}}^{\ast
})\rightarrow 0.  \label{9.1}
\end{equation}

Next, we show that $\hat{\pi}^{\ast \ast }$ is an asymptotic upper bound of $%
(1-\delta )\pi $ for any $\delta >0$. Let $\mathsf{FP}_{w^{\prime
}}(t)=\sum\nolimits_{j\in I_{0}}1_{\{|w_{j}^{\prime }|>t\}}$ and rewrite
\begin{equation*}
\bar{\Phi}_{p}(t)={\frac{\pi }{s}}\sum_{j\in I_{1}}1_{\{|\mu
	_{j}+w_{j}^{\prime }|>t\}}+{\frac{1-\pi }{p_{0}}}\mathsf{FP}_{w^{\prime
}}(t).
\end{equation*}%
Since $\hat{\pi}^{\ast \ast }>\bar{\Phi}_{p}(t)-2\bar{\Phi}(t)-c_{p}^{\ast }%
\bar{\sigma}\left( t\right) $ for any $t\in \mathbb{T}$, then
\begin{eqnarray}
\frac{\hat{\pi}^{\ast \ast }}{\pi }-1 &>&\left( {s}^{-1}\sum\nolimits_{j\in
	I_{1}}1_{\{|\mu _{j}+w_{j}^{\prime }|>t\}}-1\right) -2\bar{\Phi}(t)
\label{eqe8} \\
&&+\frac{1-\pi }{\pi }\left( p_{0}^{-1}\mathsf{FP}_{w^{\prime }}(t)-2\bar{%
	\Phi}(t)\right) -\frac{1}{\pi }c_{p}^{\ast }\bar{\sigma}(t)  \notag
\end{eqnarray}%
for any any $t\in \mathbb{T}$. Now set $t$ in the inequality (\ref{eqe8}) to
be
\begin{equation} \label{def:t_tau}
t_{\tau }=\sqrt{2\tau \log p}\text{ \ with }\tau =\gamma ^{\ast }+c/2,
\end{equation}%
where $\gamma ^{\ast }$ is defined in (\ref{def:gamma*}). We will show that
each term on the RHS of (\ref{eqe8}) is $o_{P}\left( 1\right) $.

Firstly, $c_{p}^{\ast }=O\left( \left( s_{max }/n\right) ^{1/4} \log p%
\right) $ set in \autoref{lemma:boundingSeq} implies the last term at $t_\tau$
\begin{equation*}
\pi ^{-1}c_{p}^{\ast }\bar{\sigma}(t_{\tau })=O(p^{\eta -\tau /2}(s_{max
}/n)^{1/4} \log p)=o(1).
\end{equation*}%

The second term $2\bar{\Phi}(t_{\tau })=O(p^{-\tau }/\sqrt{\log p})=o(1)$.

Consider the third term at $t_\tau$. Similar arguments for (\ref{eq:FP(xi)}) can be applied to show  $s^{-1}|\mathsf{FP}_{w'}(t_{\tau })-2p_0\bar{\Phi}(t_{\tau })| = o_P(1)$. Then
\begin{equation*}
\frac{1-\pi }{\pi }\left| p_{0}^{-1}\mathsf{FP}_{w^{\prime }}(t)-2\bar{\Phi}(t)\right| \leq Cs^{-1}|\mathsf{FP}_{w'}(t_{\tau })-2p_0\bar{\Phi}(t_{\tau })| = o_P(1) .
\end{equation*}

For the first term of (\ref{eqe8}), let $A_{1}\left( t\right) =$ ${s}^{-1}\sum\nolimits_{j\in I_{1}}1_{\{|\mu
	_{j}+w_{j}^{\prime }|\leq t\}}$. Then the first term is ${s}%
^{-1}\sum\nolimits_{j\in I_{1}}1_{\{|\mu _{j}+w_{j}^{\prime }|>t_{\tau
	}\}}-1=A_{1}\left( t_\tau\right) $. The following lemma shows $A_{1}\left(
t_{\tau }\right) =o_{P}\left( 1\right) $, and its proof is provided in %
\autoref{sec:proof_lemma_FNP(tau)}.

\begin{lemma}
	\label{lemma:FNP(tau)} Let $A_{1}\left( t\right) =$ ${s}^{-1}\sum%
	\nolimits_{j\in I_{1}}1_{\{|\mu _{j}+w_{j}^{\prime }|\leq t\}}$ and assume $%
	\mu _{min}\geq \sqrt{2(\gamma ^{\ast }+c)\log p}$. Then $A_{1}\left( t_{\tau
	}\right) =o_{P}\left( 1\right) $ for $t_{\tau }$ in (\ref{def:t_tau}).
\end{lemma}

Thus, we have shown
\begin{equation}
P(\hat{\pi}^{\ast \ast }/\pi -1<-\delta )\rightarrow 0  \label{9.2}
\end{equation}%
for any $\delta >0$. Consequently, (\ref{9}) follows from (\ref{9.1}) and (%
\ref{9.2}).

\subsection{Proof of \autoref{thm:screening}}

Recall $\mathsf{FNP}\left( t\right) =1-s^{-1}R\left( t\right)
+s^{-1}FP\left( -t\right) $ and $\widehat{\mathsf{FNP}}\left( t\right) =1-%
\hat{s}^{-1}R\left( t\right) +2\hat{s}^{-1}p\Phi \left( -t\right) $ for $%
t\geq 0$.
Recall the definition of $t^*(\epsilon)$ and simplify the notation by $t^* = t^*(\epsilon)$. We have the following \autoref{lemma:FNPdiff}, whose proof is provided in \autoref{sec:proof_lemma_FNPdiff}.
\begin{lemma}
	\label{lemma:FNPdiff}
	Assume $\mu _{min}\geq \sqrt{2(\gamma ^{\ast
		}+c)\log p}$. If $t^*$ satisfies $P\left( t^{\ast } \geq
	t_{\tau }\right) \rightarrow 1$ for $t_{\tau }$ in (\ref{def:t_tau}), then%
	\begin{equation}
	|\widehat{\mathsf{FNP}}(t^{\ast })-\mathsf{FNP}(t^{\ast })|=o_{P}(1).
	\label{eqe2.4}
	\end{equation}
\end{lemma}

Now we aim to show $P\left( t^{\ast }\geq
t_{\tau }\right) \rightarrow 1$.
The proof of the following \autoref{lemma:FNP(epsilon)} is presented in \autoref{sec:proof_lemma_FNP(epsilon)}.

\begin{lemma}
	\label{lemma:FNP(epsilon)}Assume $\mu _{min}\geq \sqrt{2(\gamma ^{\ast
		}+c)\log p}$. Then, for $t_{\tau }$ in (\ref{def:t_tau}),
	\begin{equation}
	\mathsf{FNP}(t_{\tau })=o_{P}\left( 1\right) .  \label{eqe2.2}
	\end{equation}
\end{lemma}

Note that a special case of (\ref{eqe2.4}) is $ |\widehat{\mathsf{FNP}}(t_\tau)-\mathsf{FNP}(t_\tau)|=o_{P}(1)$, which holds when $t^{\ast}$ is set to be $t_{\tau}$. Then
(\ref{eqe2.2}) implies $\widehat{\mathsf{FNP}}(t_{\tau })=o_{P}\left( 1\right)$,  and $P\left( t^{\ast } \geq t_{\tau }\right) \rightarrow 1$ follows from the definition of $t^{\ast }$.

On the other hand, the definition of $t^*$ implies $\widehat{\mathsf{FNP}}(t^{\ast
})<\epsilon $ almost surely, then \autoref{lemma:FNPdiff} implies
$P(\mathsf{FNP}(t^{\ast
})<\epsilon )\rightarrow 1$ as stated in (\ref{eq:screen}).

Next, we show (\ref{eq:screen1}). Denote
\[
\widehat{\mathsf{FNP}}_{\hat \pi} (t) = 1-\frac{R\left( t\right) -2p\Phi\left( -t\right) }{\hat{\pi} p} \qquad \text{and} \qquad \widehat{\mathsf{FNP}}_{\hat \pi^*} (t) = 1-\frac{R\left( t\right) -2p\Phi\left( -t\right) }{\hat{\pi}^* p}.
\]
By the definition of $\hat{\pi}$ and $\hat \pi^*$ and $c_p=c_p^*$, it is easy to see that $\hat \pi \ge \hat  \pi^*$ and, consequently,
\[
\widehat{\mathsf{FNP}}_{\hat \pi} (t) \ge \widehat{\mathsf{FNP}}_{\hat \pi^*} (t)
\]
for any $t>0$. Denote
\[
t^{\ast }_{\hat \pi}=\sup \left\{ t:\widehat{\mathsf{FNP}}_{\hat \pi}(t)\leq \epsilon \right\} \qquad \text{and} \qquad t^{\ast }_{\hat \pi^*}=\sup \left\{ t:\widehat{\mathsf{FNP}}_{\hat \pi^*}(t)\leq \epsilon \right\}.
\]
Then $t^*_{\hat \pi} \le t^*_{\hat \pi^*}$ almost surely.

Recall $\mathsf{FNP} (t^*_{\hat \pi^*}) < \epsilon$ with probability tending to 1 as stated in (\ref{eq:screen}) and the fact that $\mathsf{FNP} (t)$ is non-decreasing in $t$, then $\mathsf{FNP} (t^*_{\hat \pi}) \le \mathsf{FNP} (t^*_{\hat \pi^*}) < \epsilon$ with probability tending to 1. Therefore (\ref{eq:screen1}) holds.

\appendix

\section{Appendix}


The notations we will use throughout the appendices are collected as
follows. For a matrix $\mathbf{M}$, the $q$-norm $\left\Vert \mathbf{M}%
\right\Vert _{q}=\left( \sum_{i,j}\left\vert \mathbf{M}_{ij}\right\vert
^{q}\right) ^{1/q}$ for $q>0$, $\infty $-norm $\left\Vert \mathbf{M}%
\right\Vert _{\infty }=\max_{i,j}\left\vert \mathbf{M}_{ij}\right\vert $,
and $\left\Vert \mathbf{M}\right\Vert _{1,\infty }$ the maximum of the $1$%
-norm of each row of $\mathbf{M}$. If $\mathbf{M}$ is symmetric, $\sigma
_{i}(\mathbf{M})$ denotes the $i$th largest eigenvalue of $\mathbf{M}$.

\subsection{Debiased Lasso} \label{sec:extraDLasso}

The matrix $\hat{\boldsymbol{\Theta }}\in \mathbb{R}^{p\times p}$ appearing
in the debiased Lasso estimator for $\boldsymbol{\beta }$ in the main
text is obtained as follows. Let $\mathbf{X}_{-j}$ denote the matrix
obtained by removing the $j$th column of $\mathbf{X}$. For each $j=1,\ldots
,p$, let%
\begin{equation}
\hat{\boldsymbol{\gamma }}_{j}=\argmin_{\boldsymbol{\gamma }\in \mathbb{R}%
	^{p-1}}\left( n^{-1}\left\Vert \mathbf{x}_{j}-\mathbf{X}_{-j}\boldsymbol{%
	\gamma }\right\Vert _{2}^{2}+2\lambda _{j}\left\Vert \boldsymbol{\gamma }%
\right\Vert _{1}\right)  \label{eqNodewiseReg}
\end{equation}%
with components $\hat{\gamma}_{j,k},k=1,\ldots ,p$ and $k\neq j$, and define%
\begin{equation*}
\hat{\tau}_{j}^{2}=n^{-1}\left\Vert \mathbf{x}_{j}-\mathbf{X}_{-j}\hat{%
	\boldsymbol{\gamma }}_{j}\right\Vert _{2}^{2}+2\lambda _{j}\left\Vert \hat{%
	\boldsymbol{\gamma }}_{j}\right\Vert _{1}.
\end{equation*}%
Then
\begin{equation*}
\hat{\boldsymbol{\Theta }}=\text{diag}\left( \hat{\tau}_{1}^{-2},\cdots ,%
\hat{\tau}_{p}^{-2}\right) \left(
\begin{array}{cccc}
1 & -\hat{\gamma}_{1,2} & \cdots & -\hat{\gamma}_{1,p} \\
-\hat{\gamma}_{2,1} & 1 & \cdots & -\hat{\gamma}_{2,p} \\
\vdots & \vdots & \vdots & \vdots \\
-\hat{\gamma}_{p,1} & -\hat{\gamma}_{p,2} & \cdots & 1%
\end{array}%
\right) .
\end{equation*}

Recall $\sqrt{n}(\hat{\mathbf{b}}-\boldsymbol{\beta })=\mathbf{w}-%
\boldsymbol{\delta }$, where $\mathbf{w}\sim \mathcal{N}_{p}(0,\sigma ^{2}%
\hat{\boldsymbol{\Omega }})$ conditional on $\mathbf{X}$. To quantify the
magnitude of $\boldsymbol{\delta }$, we adopt and rephrase Theorem 3.13 of
\cite{JM2018} for unknown $\boldsymbol{\Sigma }$ as follows. Let $%
\boldsymbol{\Theta }=\boldsymbol{\Sigma }^{-1}$, $s_{j}=\left\vert \left\{
k\neq j:\boldsymbol{\Theta }_{jk}\neq 0\right\} \right\vert $ and $%
s_{max}=\max_{1\leq j\leq p}s_{j}$.

\begin{description}
	\item[A1)] Gaussian random design: the rows of $\mathbf{X}$ are i.i.d. $%
	\mathcal{N}_{p}\left( 0,\boldsymbol{\Sigma }\right) $ for which$\boldsymbol{%
		\ \Sigma }$ satisfies:
	
	\begin{description}
		\item[A1a)] $\max_{1\leq j\leq p}\boldsymbol{\Sigma}_{jj}\leq1$.
		
		\item[A1b)] $0<C_{min}\leq\sigma_{1}\left( \boldsymbol{\Sigma}\right)
		\leq\sigma_{p}\left( \boldsymbol{\Sigma}\right) \leq C_{max}<\infty$ for
		constants $C_{min}$ and $C_{max}.$
		
		\item[A1c)] $\rho \left( \Sigma ,C_{0}s\right) \leq \rho $ for some constant
		$\rho >0$, where $C_{0}=32C_{max }C_{min }^{-1}+1$,
		\begin{equation*}
		\rho \left( \mathbf{A},k\right) =\max_{T\subseteq \left[ p\right]
			,\left\vert T\right\vert \leq k}\big\Vert\left( \mathbf{A}_{T,T}\right)
		^{-1} \big\Vert_{1,\infty }
		\end{equation*}
		for a square matrix $\mathbf{A}$, $\left[ p\right] =\left\{ 1,...,p\right\} $
		, $\mathbf{A}_{T,T}$ is a sub-matrix formed by taking entries of $\mathbf{A}$
		whose row and column indices respectively form the same subset $T$.
	\end{description}
	
	\item[A2)] Tuning parameters: for the Lasso in \eqref{eqLasso}, $\lambda
	=8\sigma \sqrt{n^{-1}\log p}$;
	for nodewise regression in \eqref{eqNodewiseReg}, $\lambda _{j}=\tilde{\kappa%
	}\sqrt{n^{-1}\log p},j=1,\ldots ,p$ for a suitably large universal constant $%
	\tilde{\kappa}$.
	
	\item[A3)] Sparsities of $\boldsymbol{\beta }$ and $\boldsymbol{\Theta }$: $%
	\max \{s,s_{max }\}=o(n/\log p)$, $\min \{s_{max },s\}=o(\sqrt{n}/\log p)$
	and $s=o\left( n/(\log p)^{2}\right) $.
	
\end{description}

\begin{lemma}
	\label{lmJM16} Assume A1) and A2). Then there exist positive constants $c$
	and $c^{\prime }$ depending only on $C_{min }$, $C_{max }$ and $\tilde{\kappa%
	}$ such that, for $\max \{s,s_{max }\}<cn/\log p$, the probability that
	\begin{equation}
	\left\Vert \boldsymbol{\delta }\right\Vert _{\infty }\leq c^{\prime }\rho
	\sigma \sqrt{\frac{s}{n}}\log p+c^{\prime }\sigma \min \left\{ s,s_{max
	}\right\} \frac{\log p}{\sqrt{n}}  \label{eqJM16}
	\end{equation}
	is at least $1-2pe^{-16^{-1}ns^{-1}C_{min }}-pe^{-cn}-6p^{-2}$. Further, assume A3), then  $\Vert \boldsymbol{\delta }\Vert
	_{\infty }=o_{P}(1)$.
\end{lemma}

\noindent  Note that the above result relaxed the
ultra-sparse condition $s = o(\sqrt{n}/\log p)$ in \cite{vandegeer2014} to  $s=o\left( n/(\log p)^{2}\right)$ as shown in A3).


Recall the standardized debiased Lasso estimate $z_{j}=\sqrt{n}\hat{b}_{j}\sigma
^{-1}\hat{\boldsymbol{\Omega }}_{jj}^{-1/2}$ for $1\leq j\leq p$. Namely, $%
z_{j}=\mu _{j}+w_{j}^{\prime }-\delta _{j}^{\prime }$, $w_{j}^{\prime }=%
\frac{w_{j}}{\sigma \sqrt{\boldsymbol{\hat{\Omega}}_{jj}}}\sim \mathcal{N}%
(0,1)$, $\delta _{j}^{\prime }=\frac{\delta _{j}}{\sigma \sqrt{\boldsymbol{%
			\hat{\Omega}}_{jj}}}$, $\mu _{j}=\frac{\sqrt{n}\beta _{j}}{\sigma \sqrt{%
		\boldsymbol{\hat{\Omega}}_{jj}}}$ for each $j$. The $\sigma \boldsymbol{\hat{%
		\Omega}}_{jj}^{1/2}$'s \ are refereed to as standardizers. Let $\boldsymbol{%
	\delta }^{\prime }=\left( \delta _{1}^{\prime },\ldots ,\delta _{p}^{\prime
}\right) ^{T}$. We quote from \cite{Jeng:2018} some results on $\sigma
\boldsymbol{\hat{\Omega}}_{jj}^{1/2}$ for $1\leq j\leq p$ and the $%
\left\Vert \cdot \right\Vert _{1}$-norms of the covariance matrices for $%
\mathbf{w}=\left( w_{1},\ldots ,w_{p}\right) ^{T}$ and $\mathbf{w}^{\prime
}=\left( w_{1}^{\prime },\ldots ,w_{p}^{\prime }\right) ^{T}$.


\begin{lemma} \label{lmSmallop}
	Assume A2) and $s_{max }=o\left( n/\log p\right) $. Then $%
	\Vert \boldsymbol{\hat{\Omega}}-\boldsymbol{\Sigma }^{-1}\Vert _{\infty
	}=o_{P}\left( 1\right) $. If further A1b) holds, then $\Vert \hat{
		\boldsymbol{\Theta }}\hat{\boldsymbol{\Sigma }}-\mathbf{I}\Vert _{\infty
	}=O_{P}(\lambda _{1})$, both $\min_{1\leq j\leq p}\hat{\boldsymbol{\Omega }}%
	_{jj}$ and $\max_{1\leq j\leq p}\hat{\boldsymbol{\Omega }}_{jj}$ are
	uniformly bounded (in $p$) away from $0$ and $\infty $ with probability
	tending to $1$, and $\left\Vert \boldsymbol{\delta }^{\prime }\right\Vert
	_{\infty }\leq (\sigma \sqrt{C_{min }})^{-1}\left\Vert \boldsymbol{\delta }%
	\right\Vert _{\infty }$ with probability tending to $1$.
\end{lemma}

\begin{lemma}
	\label{lemma:K} Let $\hat{\mathbf{K}}$ be the correlation matrix of $\mathbf{%
		\ w}$. Assume A1) and A2). Then
	\begin{equation}
	p^{-2}\Vert \sigma ^{2}\hat{\boldsymbol{\Omega }}\Vert
	_{1}=O_{P}\left(\lambda _{1}\sqrt{s_{max }}\right) \quad \text{and}\quad
	\Vert \hat{\mathbf{\ K}}\Vert _{1}=O(\sigma ^{2}\Vert \hat{\boldsymbol{%
			\Omega }}\Vert _{1}).  \label{eqCorMatL1}
	\end{equation}
\end{lemma}

\subsection{Hermite polynomials and Mehler expansion}

\label{secHermitePoly}

The following is quoted from \cite{Jeng:2018}. Let $\phi \left( x\right)
=\left( 2\pi \right) ^{-1/2}\exp \left( -x^{2}/2\right) $ and
\begin{equation*}
f_{\rho }\left( x,y\right) =\frac{1}{2\pi \sqrt{1-\rho ^{2}}}\exp \left( -%
\frac{x^{2}+y^{2}-2\rho xy}{2\left( 1-\rho ^{2}\right) }\right)
\end{equation*}%
for $\rho \in \left( -1,1\right) $. For a nonnegative integer $k$, let $%
H_{k}\left( x\right) =\left( -1\right) ^{k}\frac{1}{\phi \left( x\right) }%
\frac{d^{k}}{dx^{k}}\phi \left( x\right) $ be the $k$th Hermite polynomial;
see \cite{Feller:1971B} for such a definition. Then Mehler's expansion %
\citep{Mehler:1866} gives%
\begin{equation}
f_{\rho }\left( x,y\right) =\left( 1+\sum\nolimits_{k=1}^{\infty }\frac{\rho
	^{k}}{k!}H_{k}\left( x\right) H_{k}\left( y\right) \right) \phi \left(
x\right) \phi \left( y\right) .  \label{eq:Mehler}
\end{equation}%
Further, Lemma 3.1 of \cite{chen2016} asserts%
\begin{equation}
\left\vert e^{-y^{2}/2}H_{k}\left( y\right) \right\vert \leq C_{0}\sqrt{k!}%
k^{-1/12}e^{-y^{2}/4}\text{ \ for any\ }y\in \mathbb{R}
\label{eqBoundHermite}
\end{equation}%
for some constant $C_{0}>0$.

	\label{AppMain}
	
	\subsection{Proof of \autoref{lemma:zw-difference}} \label{sec:proof_zw-difference}
	By assumption A3), $d_{p}=o(1)$, $s\ll n/\log p$, and $%
	n\gg \log p$. Then \autoref{lmJM16} implies
	\begin{equation*}
	P(\Vert \boldsymbol{\delta }\Vert _{\infty }\geq d_{p})\leq 2pe^{-c_{\ast
		}n/s}+pe^{-Cn}+6p^{-2}\leq Cp^{-2}
	\end{equation*}%
	where $c_{\ast }=C_{\min }/16$. By \autoref{lmSmallop} and the definition of
	$d_{p}$, we have $P(\Vert \boldsymbol{\delta }^{\prime }\Vert _{\infty }\geq
	d_{p})\leq Cp^{-2}$.
	
	Now consider $ P(|w'_{j} - \delta'_j|>t_{\xi})$, which is bounded as follows.
	\[
	P(|w'_{j}| >t_{\xi} + |\delta'_j|) \le P(|w'_{j} - \delta'_j|>t_{\xi}) \le  P(|w'_{j}| >t_{\xi} - |\delta'_j|).
	\]
	The rightmost term
	\begin{eqnarray*}
		P(|w'_{j}| >t_{\xi} - |\delta'_j|) & \le & P(|w'_{j}| >t_{\xi} - |\delta'_j|, \Vert\ms{\delta}' \Vert_\infty \le d_p) + P(\Vert\ms{\delta}' \Vert_\infty > d_p)
		\\
		& \le &  P(|w'_{j}| >t_{\xi} - d_p) + Cp^{-2}.
	\end{eqnarray*}
	On the other hand, the leftmost term
	\begin{eqnarray*}
		P(|w'_{j}| >t_{\xi} + |\delta'_j|) & \ge & P(|w'_{j}| >t_{\xi} + |\delta'_j|, \Vert\ms{\delta}' \Vert_\infty \le d_p)	\\
		& \ge &  P(|w'_{j}| >t_{\xi} + d_p, \Vert\ms{\delta}' \Vert_\infty \le d_p) \\	
		& = &  P(|w'_{j}| >t_{\xi} + d_p) - P(|w'_{j}| >t_{\xi} + d_p, \Vert\ms{\delta}' \Vert_\infty > d_p) \\
		& \ge &  P(|w'_{j}| >t_{\xi} + d_p) - P(\Vert\ms{\delta}' \Vert_\infty > d_p)  \\
		& \ge & P(|w'_{j}| >t_{\xi} + d_p) - Cp^{-2}.
	\end{eqnarray*}
	Summing up the above gives
	\[
	\left\vert P(|w'_{j} - \delta'_j|>t_{\xi}) - P(|w'_{j}|>t_{\xi })\right\vert \le C \phi(t_\xi) d_p + C p^{-2},
	\]
	and the claim in \autoref{lemma:zw-difference} follows.

	\subsection{Proof of \autoref{lemma:Var(FP)}}
	
	\label{sec:proof_lemma_Var(FP)}
	
	For $i\neq j$, let $\rho _{ij}$ be the correlation
	between $w_{i}^{\prime }$ and $w_{j}^{\prime }$ and $C_{ij,\xi }=\mathsf{Cov}%
	\left( 1_{\left\{ |w_{i}^{\prime }|\leq t_{\xi }\right\} },1_{\left\{
		|w_{j}^{\prime }|\leq t_{\xi }\right\} }\right) $. Then, by \autoref%
	{lmSmallop}, $\rho _{ij}$ is also the correlation between $w_i$ and $w_j$.
	Further,
	\begin{equation}
	\mathsf{Var}\left(
	\sum\nolimits_{j=1}^{p}1_{\{|w_{j}^{\prime }|>t_{\xi }\}}\right) \leq \sum_{j=1}^p
	\mathsf{Var}\left( 1_{\left\{ \left\vert w_{j}^{\prime }\right\vert \leq
		t_{\xi }\right\} }\right) +\sum_{i\neq j}C_{ij,\xi }.  \label{0}
	\end{equation}%
	By Mill's ratio,
	\begin{equation}
	\sum_{j=1}^p\mathsf{Var}\left( 1_{\left\{ \left\vert w_{j}^{\prime
		}\right\vert \leq t_{\xi }\right\} }\right) \leq 2p\Phi \left( -t_{\xi
	}\right) \left( 1-2\Phi \left( -t_{\xi }\right) \right) =O(p^{1-\xi }).
	\label{1.1}
	\end{equation}%
	It is left to bound $\sum_{i\neq j}C_{ij,\xi }$ in (\ref{0}).
	
	Define $c_{1,\xi }=-t_{\xi }$ and $c_{2,\xi }=t_{\xi }$. Fix a pair of $%
	(i,j) $ such that $i\neq j$ and $|\rho _{ij}|\neq 1$. Now we will use the
	results in \autoref{secHermitePoly}. Since $C_{ij,\xi }$ is finite and the
	series in Mehler's expansion in (\ref{eq:Mehler}) as a trivariate function
	of $\left( x,y,\rho \right) $ is uniformly convergent on each compact set of
	$\mathbb{R}\times \mathbb{R}\times \left( -1,1\right) $ as justified by \cite%
	{Watson:1933}, we can interchange the order of the summation and integration
	and obtain%
	\begin{eqnarray*}
		C_{ij,\xi } &=&\int_{c_{1,\xi }}^{c_{2,\xi }}\int_{c_{1,\xi }}^{c_{2,\xi
			}}f_{\rho _{ij}}\left( x,y\right) dxdy-\int_{c_{1,\xi }}^{c_{2,\xi }}\phi
			(x)dx\int_{c_{1,\xi }}^{c_{2,\xi }}\phi (y)dy \\
			&=&\sum_{k=1}^{\infty }\frac{\rho _{ij}^{k}}{k!}\int_{c_{1,\xi }}^{c_{2,\xi
				}}H_{k}(x)\phi (x)dx\int_{c_{1,\xi }}^{c_{2,\xi }}H_{k}(y)\phi (y)dy.
			\end{eqnarray*}%
			Since $H_{k-1}\left( x\right) \phi \left( x\right) =\int_{-\infty
			}^{x}H_{k}\left( y\right) \phi \left( y\right) dy$ for $x\in \mathbb{R}$,
			then
			\begin{eqnarray*}
				C_{ij,\xi } &=&\sum_{k=1}^{\infty }\frac{\rho _{ij}^{k}}{k!}\left[
				H_{k-1}(c_{2,\xi })\phi (c_{2,\xi })-H_{k-1}(c_{1,\xi })\phi (c_{1,\xi })%
				\right] ^{2} \\
				&\leq &2\sum_{k=1}^{\infty }\frac{\left\vert \rho _{ij}\right\vert ^{k}}{k!}%
				\left\{ \left[ H_{k-1}\left( c_{2,\xi }\right) \phi \left( c_{2,\xi }\right) %
				\right] ^{2}+\left[ H_{k-1}\left( c_{1,\xi }\right) \phi \left( c_{1,\xi
				}\right) \right] ^{2}\right\} .
			\end{eqnarray*}%
			Inequality (\ref{eqBoundHermite}) implies, for some finite constant $C_{0}>0$
			,
			\begin{equation*}
			\left[ H_{k-1}\left( c_{2,\xi }\right) \phi \left( c_{2,\xi }\right) \right]
			^{2}+\left[ H_{k-1}\left( c_{1,\xi }\right) \phi \left( c_{1,\xi }\right) %
			\right] ^{2}\leq C_{0}^{2}(k-1)!(k-1)^{-1/6}e^{-t_{\xi }^{2}/2}.
			\end{equation*}%
			Therefore,
			\begin{eqnarray}
			\left\vert \sum\nolimits_{i\neq j}C_{ij,\xi }\right\vert &\leq &C\sum_{1\leq
				i<j\leq p}|\rho _{ij}|\sum_{k=1}^{\infty }k^{-7/6}\left\vert \rho
			_{ij}\right\vert ^{k-1}e^{-t_{\xi }^{2}/2}  \notag \\
			&\leq &Cp^{-\xi }\sum_{1\leq i<j\leq p}|\rho _{ij}|=O(p^{-\xi }\Vert \hat{%
				\mathbf{K}}\Vert _{1}).  \label{1.2}
			\end{eqnarray}%
			Combining (\ref{0}) with (\ref{1.1}) and (\ref{1.2}) gives
			\begin{equation*}
			\mathsf{Var}\left(
			\sum\nolimits_{j=1}^{p}1_{\{|w_{j}^{\prime }|>t_{\xi }\}}\right) =O(p^{1-\xi })+O(p^{-\xi
			}\Vert \hat{\mathbf{K}}\Vert _{1})=O(p^{1-\xi })+O(p^{2-\xi }\lambda _{1}%
			\sqrt{s_{max }}),
			\end{equation*}%
			where the last inequality follows from \autoref{lemma:K}, i.e., $\Vert
			\hat{\mathbf{K}}\Vert _{1}=O_{P}(p^{2}\lambda _{1}\sqrt{s_{max }})$.
			

\subsection{Proof of \autoref{lemma:boundingSeq}} \label{sec:proof_lemma_boundingSeq}

Recall $\bar{\sigma}\left( t\right) =\sqrt{2\bar{\Phi}(t)\left( 1-2\bar{\Phi}%
	(t)\right) }$ and $\mathsf{H}(t) =\left( \bar{\sigma}\left( t\right) \right)
^{-1}\left( p^{-1}\sum\nolimits_{j=1}^{p}1_{\{|w_{j}^{\prime }|>t\}}-2\bar{%
	\Phi}(t)\right) $. Then $E(\mathsf{H}(t) )=0$ since $w_{j}^{\prime }\sim
\mathcal{N}_{1}(0,1)$ for all $j$.

For any $t_{\xi }=\sqrt{2\xi \log p}$ such that $\lim_{p\rightarrow \infty }t_{\xi }=\infty $, \autoref{lemma:Var(FP)} implies%
\begin{eqnarray}
\mathsf{Var}\left( \mathsf{HC}(t_{\xi })\right) &=&p^{-2}\bar{\sigma}%
_{p}^{-2}\left( t_{\xi }\right) \mathsf{Var}\left(
\sum\nolimits_{j=1}^{p}1_{\{|w_{j}^{\prime }|>t_{\xi }\}}\right)  \notag \\
&\leq &Cp^{\xi -2}\sqrt{\log p}\left( p^{1-\xi }+p^{2-\xi }\sqrt{\log p}%
\sqrt{s_{max }/n}\right)  \notag \\
&=&O\left(\sqrt{s_{max }/n} \log p \right) .  \label{eqa7}
\end{eqnarray}

Let
\begin{equation} \label{def:T}
\mathbb{T}=\left[ \sqrt{\tau _{0}\log p},\sqrt{\tau _{1}\log p}\right]
\cap \mathbb{N}
\end{equation}
for which $0<\tau _{0}<\tau _{1}$. So, each $t\in \mathbb{T}
$ can be written as $t=t_{\xi }=\sqrt{2\xi \log p}$ for some $\xi >0$ and $%
\lim_{p\rightarrow \infty }t_{\xi }=\infty $. Recall $V_{p}^{\ast }=\max
\left\{ \mathsf{H}(t)  :t\in \mathbb{T}\right\} $. Therefore, (%
\ref{eqa7}) implies
\begin{eqnarray*}
	P(V_{p}^{\ast }>c_{p}^{\ast }) &\leq &C\left( c_{p}^{\ast }\right) ^{-2}%
	\sqrt{\log p}\max_{t\in \mathbb{T}}\mathsf{Var}\left( \mathsf{H}(t) \right)
	\\
	&\leq &C\left( c_{p}^{\ast }\right) ^{-2} \cdot (\log p)^{3/2} \cdot \sqrt{s_{max }/n}.
\end{eqnarray*}%
However, $c_{p}^{\ast }=O\left( \left( s_{max }/n\right) ^{1/4} \log p%
\right) $. Thus, $P(V_{p}^{\ast }>c_{p}^{\ast })=o\left( 1\right) $ as
desired. \newline			
			
			\subsection{Proof of \autoref{lmLast}} \label{sec:proof_empiricalDiff}
			Since $\max_{t\in
				\mathbb{T}}\left( 1-2\bar{\Phi}(t)\right) \geq 4^{-1}$ for all $p$
			sufficiently large. It suffices to show
			\begin{equation}
			\pi ^{-1}\left\vert F_{p}(t)-\bar{\Phi}%
			_{p}(t)\right\vert =o_{P}(1)  \label{eqe12}
			\end{equation}%
			for $t = t_{\xi }=\sqrt{2\xi \log p}$ with $\xi >\gamma^* $. Perform the decomposition
			\[
			\left\vert F_{p}(t)-\bar{\Phi}_{p}(t)\right\vert \le \left\vert F_{p}(t)-E(F_p(t))\right\vert + \left\vert \bar{\Phi}_{p}(t) - E(\bar{\Phi}_{p}(t))\right\vert + \left\vert E(F_p(t))- E(\bar{\Phi}_{p}(t))\right\vert.
			\]
			Similar arguments for (\ref{eq:1}) can be applied to show
			\[
			\pi^{-1}\left\vert F_{p}(t)-E(F_p(t))\right\vert = o_p(1) = \pi^{-1}\left\vert \bar{\Phi}_{p}(t) - E(\bar{\Phi}_{p}(t))\right\vert,
			\]
			and similar arguments for (\ref{eq:2}) imply
			\[
			\pi^{-1}\left\vert E(F_p(t))- E(\bar{\Phi}_{p}(t))\right\vert = o(1).
			\]
			Summing up the above gives (\ref{eqe12}).

			\subsection{Proof of \autoref{lemma:FNP(tau)}}
			
			\label{sec:proof_lemma_FNP(tau)}
			
			We will show $A_{1}\left(
			t_{\tau }\right) =o_{P}\left( 1\right) $. Fix a constant $a>0$,
			\begin{equation}
			P\left( A_{1}\left( t_{\tau }\right) >a\right) \leq \frac{1}{as}\sum_{j\in
				I_{1}}P\left( |\mu _{j}+w_{j}^{\prime }|\leq t_{\tau }\right) \leq \frac{1}{a%
			}\max_{j\in I_{1}}P\left( |\mu _{j}+w_{j}^{\prime }|\leq t_{\tau }\right)
			\label{eqe15}
			\end{equation}%
			and for each $j\in I_{1}$%
			\begin{equation}
			P\left( |\mu _{j}+w_{j}^{\prime }|\leq t_{\tau }\right) =1-\bar{\Phi}\left(
			t_{\tau }-\mu _{j}\right) -\Phi \left( -t_{\tau }-\mu _{j}\right) ,
			\label{eqe14}
			\end{equation}%
			We only need to uniformly bound the RHS of (\ref{eqe14}).
			
			Recall
			$\mu _{min}=\min_{j\in I_{1}}\sqrt{n}\left\vert \beta _{j}\right\vert \sigma
			^{-1}\sqrt{\boldsymbol{\Sigma }_{jj}}$ and condition $\mu_{min} \ge \sqrt{2(\gamma^*+c) \log p}$. Let
			\[
			\mu _{min}=\sqrt{2r\log p},
			\]
			then  $r \ge \tau +c/2$. Further, by \autoref{lmSmallop}, the ratio $\mu _{min} / \min_{j\in
				I_{1}}\left\vert \mu _{j}\right\vert $ is uniformly bounded (in $p$) away
			from $0$ and $\infty $. Then, two cases happen for each $j\in I_{1}$: (i) both $t_{\tau
			}-\mu _{j}\rightarrow -\infty $ and $-t_{\tau }-\mu _{j}\rightarrow -\infty $
			when $\mu _{j}>0$; (b) both $t_{\tau }-\mu _{j}\rightarrow +\infty $ and $%
			-t_{\tau }-\mu _{j}\rightarrow +\infty $ when $\mu _{j}<0$. However, in
			either case,%
			\begin{equation*}
			\min_{j\in I_{1}}\min \left\{ \left\vert t_{\tau }-\mu _{j}\right\vert
			,\left\vert -t_{\tau }-\mu _{j}\right\vert \right\} \geq \sqrt{2\tilde{c}%
				\log p},
			\end{equation*}%
			where $\tilde{c}=2^{-1}\left( \sqrt{2\tau +c}-\sqrt{2\tau }\right) ^{2}>0$.
			Therefore,
			\begin{equation}
			\max_{j\in I_{1}}P\left( |\mu _{j}+w_{j}^{\prime }|\leq t_{\tau }\right)
			\leq 4\bar{\Phi}\left( \sqrt{2\tilde{c}\log p}\right) =O\left( p^{-\tilde{c}%
			}\right) .  \label{eqe16}
			\end{equation}%
			Combining (\ref{eqe16}) with (\ref{eqe15}) gives%
			\begin{equation*}
			P\left( A_{1}\left( t_{\tau }\right) >a\right) \leq a^{-1}O\left( p^{-\tilde{%
					c}}\right) =o\left( 1\right) ,
			\end{equation*}%
			which is the desired claim on $A_{1}\left( t_{\tau }\right) $.
			
			\subsection{Proof of \autoref{lemma:FNPdiff}} \label{sec:proof_lemma_FNPdiff}
			Recall $\widehat{\mathsf{FNP}}\left( t\right) =1-\hat{s}^{-1}\left(
			R\left( t\right) -2(p-\hat s)\Phi \left( -t\right) \right) $. We only need to show
			\begin{eqnarray} \label{eq:3}
			&&\left\vert s^{-1}\left( R\left( t^{\ast }\right) -\mathsf{FP}\left(
			t^{\ast }\right) \right) -\hat{s}^{-1}\left( R\left( t^{\ast }\right)
			-2(p-\hat s)\Phi \left( -t^{\ast }\right) \right) \right\vert \nonumber \\
			&\le& \left\vert \mathsf{TP}\left( t^{\ast }\right) \left( s^{-1}-\hat{s}
			^{-1}\right) \right\vert + \left\vert \hat{s}^{-1} (\mathsf{FP}\left( t^{\ast }\right) - 2p_0\Phi \left( -t^{\ast }\right) \right\vert  \\
			& + &\left\vert 2 \hat s^{-1} \Phi(-t^*) (p-\hat s- p_0)\right\vert =o_{P}\left( 1\right). \nonumber
			\end{eqnarray}
			Since $\mu _{min}\geq \sqrt{2(\gamma ^{\ast }+c)\log p}$, then  \autoref{thm:proportion} implies
			\begin{equation}
			P\left( 1-\delta \leq \hat{s}s^{-1}\leq 1\right) \rightarrow 1  \label{eqe20}
			\end{equation}
			for any $\delta >0$. Let $\delta ^{\prime }=\frac{\delta }{1-\delta }$. Then, (\ref{eqe20}) is
			equivalent to%
			\begin{equation*}
			P\left( 0\leq \hat{s}^{-1}-s^{-1}\leq \delta ^{\prime }s^{-1}\right)
			\rightarrow 1.
			\end{equation*}
			Pick a $\delta >0$ such that $\delta <\frac{a}{1+a}$. Then $\delta <1$ and $%
			\delta ^{\prime }s^{-1}\mathsf{TP}\left( t^{\ast }\right) <a$ almost surely.
			Therefore,
			\begin{eqnarray*}
				&&P\left( \left\vert \mathsf{TP}\left( t^{\ast }\right) \left( s^{-1}-\hat{s}%
				^{-1}\right) \right\vert >a\right)  \\
				&\leq &P\left( \left\vert \mathsf{TP}\left( t^{\ast }\right) \left( s^{-1}-%
				\hat{s}^{-1}\right) \right\vert >a,\left\vert s^{-1}-\hat{s}^{-1}\right\vert
				\leq \delta ^{\prime }s^{-1}\right) +P\left( \left\vert s^{-1}-\hat{s}%
				^{-1}\right\vert >\delta ^{\prime }s^{-1}\right)  \\
				&\leq &P\left( \mathsf{TP}\left( t^{\ast }\right) \delta ^{\prime
				}s^{-1}\geq a\right) +o\left( 1\right)  \\
				&=&0+o\left( 1\right) ,
			\end{eqnarray*}%
			i.e., the first term in (\ref{eq:3}) $= o_P(1)$. The remaining two terms in (\ref{eq:3}) are also of $ o_P(1)$ by (\ref{eqe20}) and \autoref{thm:FP}. This concludes the proof.

			\subsection{Proof of \autoref{lemma:FNP(epsilon)}}
			
			\label{sec:proof_lemma_FNP(epsilon)} Recall $\mathsf{FNP}\left( t\right)
			=s^{-1}\sum_{j\in I_{1}}1_{\left\{ \left\vert \mu _{j}+w_{j}^{\prime
				}-\delta _{j}^{\prime }\right\vert \leq t\right\} }$ and $A_{1}\left(
			t\right) =$ ${s}^{-1}\sum\nolimits_{j\in I_{1}}1_{\{|\mu _{j}+w_{j}^{\prime
				}|\leq t\}}$. Since $\mu _{min}\geq \sqrt{2(\gamma ^{\ast }+c)\log p}$, \autoref{lemma:FNP(tau)}
			implies $A_{1}\left( t_{\tau }\right) =o_{P}\left( 1\right) $. Now we show $%
			\mathsf{FNP}(t_{\tau })=o_{P}\left( 1\right) $. Clearly,%
			\begin{equation*}
			P\left( \mathsf{FNP}(t_{\tau })>a\right) \leq \frac{1}{a}\max_{j\in
				I_{1}}P\left( |\mu _{j}+w_{j}^{\prime }+\delta _{j}^{\prime }|\leq t_{\tau
			}\right) .
			\end{equation*}%
			However, $\max_{1\leq i\leq p}|\delta _{j}^{\prime }|=o_{P}\left( 1\right) $
			and $w_{j}^{\prime }\sim \mathcal{N}_{1}(0,1)$ for each $j$ together imply%
			\begin{equation*}
			\max_{j\in I_{1}}\left\vert P\left( |\mu _{j}+w_{j}^{\prime }+\delta
			_{j}^{\prime }|\leq t_{\tau }\right) -P\left( |\mu _{j}+w_{j}^{\prime }|\leq
			t_{\tau }\right) \right\vert =o\left( 1\right) .
			\end{equation*}%
			Combining the above with (\ref{eqe16}) gives $\max_{j\in I_{1}}P\left( |\mu
			_{j}+w_{j}^{\prime }+\delta _{j}^{\prime }|\leq t_{\tau }\right) =o\left(
			1\right) $, and $\mathsf{FNP}(t_{\tau })=o_{P}\left( 1\right) $ holds.

			\bibliographystyle{chicago}

\end{document}